\documentclass[12pt,amssymb]{article}
\usepackage{amsfonts,latexsym}

\newcommand{\proj}{\mathbb P}
\newcommand{\reel}{\mathbb R}
\newcommand{\complex}{\mathbb C}

\newcommand{\é}{\'{e}}

\newcommand{\à}{\`{a}}
\newcommand{\ù}{\`{u}}

\newtheorem{theorem}{Th\éor\`{e}me}
\newtheorem{lem}{Lemme}

\newtheorem{coro}{Corollaire}

\begin{document}
\title{Enveloppes infé\érieures de fonctions admissibles sur l'espace projectif complexe.
 Cas dissym\étrique.}
\author {\em Adn\ène Ben Abdesselem et Brahim Dridi}
\date{\em Universit\é Paris 6, U.F.R.929, 4, place Jussieu, 75005, Paris,\\  
et L.I.M. Ecole Polytechnique de Tunisie.}

\maketitle {\em RESUME. Cet article g\én\éralise les
r\ésultats du premier auteur concernant les fonctions
admissibles sur certaines vari\ét\és de Fano
\cite{B},\cite{B1}. On \étudie ici une classe plus large de
fonctions, ces derni\ères pouvant pr\ésenter un d\éfaut
de sym\étrie, et on montre l'existence d'une fonction limite
donnant pr\écis\ément l'invariant de Tian sur les
vari\ét\és consid\ér\ées et minorant toutes les
fonctions admissibles dont le sup est \égal \`{a} z\éro.}

\maketitle {\em ABSTRACT. This paper generalizes the first
author's preceding works concerning admissible functions on
certain Fano manifolds \cite{B},  \cite{B1}. Here, we study a
larger class of functions which can be less symmetric than the
ones studied before. When the sup of these functions is null, we
prove that they admit a lower bound, giving precisely Tian
invariant on these manifolds. } 
{\section{Introduction.}\label{i}} Dans un travail
pr\éc\édent \cite{B1}, le premier auteur a prouv\é
l'existence d'une fonction minorant toutes les fonctions
admissibles \`{a} sup \égal \`{a} z\éro sur certains
projectifs \éclat\és et invariantes par un groupe
d'automorphismes bien choisi. Ce dernier s'obtient \`{a} partir de
celui de $\proj_{m}\complex$ et contient en particulier les
automorphismes \échangeant deux coordonn\ées homog\ènes
entre elles. Dans cet article les projectifs \éclat\és
consid\ér\ées ne peuvent, de par leur g\éom\étrie,
avoir un groupe d'automorphismes aussi riche. On montre alors le
r\ésultat ci-dessous expos\é sur $\proj_{m}\complex$ en
consid\érant le groupe d'automorphismes adé\équat, nous
permettant, \`{a} l'instar de ce qui a \ét\é fait dans
\cite{B1}, d'é\étendre la mé\éthode aux espaces X et Y
é\étudié\és dans  \cite{BC}.

Rappelons dans un premier temps la d\éfinition des espaces,
mé\étriques, et groupes d'automorphismes que l'on
considé\érera. On adoptera les notations utilisé\ées dans
\cite{BC}. D\ésignons par $[z_{0},z_{1},..,z_{m}]$ les
coordonn\ées homog\ènes de l'espace projectif complexe
$\proj_{m}\complex$
de dimension complexe $m\geq2$.\\
Munissons $\proj_{m}\complex$ de la m\étrique $g$ ayant pour
composantes, dans la carte $\{z_{0}\neq 0\}$,
$$g_{\lambda\overline{\mu}}=(m+1){\partial}_{\lambda
\bar{\mu}}\ln (1+ x_{1}+..+x_{m})$$ o\`{u} $x_{i}=\mid z_{i}
\mid^{2}$ et ${\partial}_{\lambda
\bar{\mu}}=\frac{\partial^{2}}{\partial
z_{\lambda}\partial\overline{z}_{\mu}}$ .

Fixons pour le reste de l'article l'entier $k\in\{1,...,m-1\}$, et
rappelons la d\éfinition des vari\ét\és $X$ et $Y$
utilis\ées dans \cite{BC}.

$X$ est l'é\éclaté\é de $\proj_{m}\complex$ au dessus du
sous-ensemble $\{[0,..,0,z_{k+1},..,z_{m}]\}$ qui s'identifie
\`{a} l'espace $\proj_{m-k-1}\complex$; $Y$ est l'\éclat\é
de $\proj_{m}\complex$ au dessus des sous-ensembles
$\{[z_{0},..,z_{k},0,..,0]\}$ et $\{[0,..,0,z_{k+1},..,z_{m}]\}$
qui s'identifient respectivement \`{a} $\proj_{k}\complex$ et
$\proj_{m-k-1}\complex$.

$X$ est alors la sous varié\été\é de $\proj_{m}\complex
\times\proj_{k}\complex$ constitu\ée des points
$$([z_{0},..,z_{k},z_{k+1},..,z_{m}],[\zeta_{0},..,\zeta_{k}])\in
\proj_{m}\complex \times\proj_{k}\complex$$ tels que les deux
vecteurs $(z_{0},..,z_{k})$ et $(\zeta_{0},..,\zeta_{k})$ de
$\complex^{k+1}$ soient colin\éaires.

On consid\ère alors les projections $\pi_{1}$ et $\pi_{2}$ de
$X$ respectivement sur $\proj_{m}\complex$ et $\proj_{k}\complex$.
En utilisant les m\étriques de Fubini-Study $g_{m}$ de
$\proj_{m}\complex$ et $g_{k}$ de $\proj_{k}\complex$, on
d\éfinit la m\étrique $\tilde{g}$ sur $X$ par
$$\tilde{g}=(m+1-k)\pi_{1}^{*}g_{m} + k\pi_{2}^{*}g_{k}.$$ Ses composantes
dans la carte dense de $X$ constitu\ée des points :
$$([1,z_{1},..,z_{m}], [1,z_{1},..,z_{k}]); (z_{1},..,z_{m}) \in
\complex ^{m},$$ sont donn\ées par
$$\tilde{g}_{\lambda\overline{\mu}}=(m+1-k){\partial}_{\lambda
\bar{\mu}}\ln (1+ x_{1}+..+x_{m}) + k{\partial}_{\lambda
\bar{\mu}}\ln (1+ x_{1}+..+x_{k})$$

De m\^{e}me, $Y$ est la sous-vari\été\é de
$\proj_{m}\complex\times\proj_{k}\complex\times\proj_{m-k-1}\complex$
constitu\ée de points
$$([z_{0},..,z_{k},z_{k+1},..,z_{m}],[\zeta_{0},..,\zeta_{k}],[\eta_{k+1},..,\eta_{m}])
\in\proj_{m}\complex\times\proj_{k}\complex\times\proj_{m-k-1}\complex$$
tels que les vecteurs $(z_{0},..,z_{k})$ et $(z_{k+1},..,z_{m})$
sont respectivement proportionnels \`{a}
$(\zeta_{0},..,\zeta_{k})$ et $(\eta_{k+1},..,\eta_{m})$. Comme
pour $X$, on consid\ère les projections $\pi_{1}$, $\pi_{2}$ et
$\pi_{3}$ de $Y$ respectivement sur $\proj_{m}\complex$,
$\proj_{k}\complex$ et $\proj_{m-k-1}\complex$. Alors
$$\hat{g}=2\pi_{1}^{*}g_{m} +
k\pi_{2}^{*}g_{k}+(m-k-1)\pi_{3}^{*}g_{m-k-1}$$ est une
m\étrique dans $Y$ qui a pour composantes
$$\hat{g}_{\lambda\overline{\mu}}=2{\partial}_{\lambda
\bar{\mu}}\ln (1+ x_{1}+..+x_{m}) + k{\partial}_{\lambda
\bar{\mu}}\ln (1+ x_{1}+..+x_{k}) + (m-k-1){\partial}_{\lambda
\bar{\mu}}\ln (x_{k+1}+..+x_{m})$$ dans la carte dense de $Y$
$$\{([1,z_{1},..,z_{m}],[1,z_{1},..,z_{k}],[z_{k+1},..,z_{m}]),(z_{1},..,z_{m})\in
\complex^{m} \mbox{ et } (z_{k+1},..,z_{m})\neq 0 \}$$

On rappelle que $\tilde{g}$ et $\hat{g}$ sont respectivement dans
la premi\ère classe de Chern de $X$ et de $Y$ (voir \cite{BC});
et par cons\équent, $X$ et $Y$ sont de Fano.

On consid\ère le groupe d'automorphismes $G$ sur
$\proj_{m}\complex$ engendr\é par les automorphismes
$\sigma_{i,j}$, $\varphi_{p,q}$ $\tau_{l,\theta}$ d\éfinies
$\forall$ $i,j \in \{0,1,..,k\}$, $p,q \in \{k+1,..,m\}$, $l
\in\{0,1,..,m\}$ et $\theta \in [0,2\pi]$ par
\\$$\sigma_{i,j}([z_{0},..,z_{i},..,z_{j},..,z_{k},..z_{m}])
=[z_{0},..,z_{j},..,z_{i},..,z_{k},..,z_{m}]$$
$$\varphi_{p,q}([z_{0},..,z_{k},..,z_{p},..,z_{q},..z_{m}])
=[z_{0},..,z_{k},..,z_{q},..,z_{p},..z_{m}]$$ et
$$\tau_{l,\theta}([z_{0},..,z_{l},..,z_{m}])=
[z_{0},..,z_{l}e^{i\theta},..,z_{m}].$$
\\Ce groupe engendre des groupes d'automorphismes
naturels de $X$ et $Y$ que l'on notera encore $G$ dans les deux
cas.

On d\éfinit sur
$\complex^{m+1}\backslash\bigcup_{p}\{z_{p}=0\}$ la fonction
 $\psi =\inf(\psi_1,\psi_2)$, o\`{u}
$$\psi_1=\ln \frac{(\mid z_{0}\mid...\mid z_{k}\mid
)^{2{(m+1)}/(k+1)}}{(\mid z_{0}\mid^{2}+...+\mid
z_{m}\mid^{2})^{(m+1)}}$$ et $$\psi_2=\ln \frac{(\mid
z_{k+1}\mid...\mid z_{m}\mid )^{2{(m+1)}/(m-k)}}{(\mid
z_{0}\mid^{2}+...+\mid z_{m}\mid^{2})^{(m+1)}}.$$\\
$\psi_1$ et $\psi_2$ sont homog\ènes de degr\é z\éro sur
$\complex^{m+1}$, chacune d'entre elles induit alors une fonction
sur $\proj_{m}\complex$. La fonction $\psi_1$ atteint son maximum
\'egal \`a $-(m+1)\ln (k+1)$ en les points
$[1,e^{i\theta_{1}},..,e^{i\theta_{k}},0,..,0]
\in\proj_{m}\complex$ et tend vers moins l'infini lorsque l'une
des coordonn\ées homog\ènes $z_{0},..,z_{k}$ tend vers
z\éro ou vers l'infini, ce qui correspond aux fronti\ères
des cartes denses d\éfinies par $\{z_{i}\neq 0, 0\leq i\leq k\}$.\\
$\psi_{2}$ atteint son maximum \'egal \`a $-(m+1)\ln (m-k)$ en les
points\\ $[0,..,0,
e^{i\theta_{k+1}},..,e^{i\theta_{m}}]\in\proj_{m}\complex$ et tend
vers moins l'infini lorsque l'une des coordonn\ées
homog\ènes $z_{k+1},..,z_{m}$ tend vers z\éro ou vers
l'infini, c'est-\`{a}-dire aux fronti\ères des cartes denses
d\éfinies par $\{z_{j}\neq 0, k+1\leq j\leq m\}$.

Pour d\écrire la fonction extr\émale $\tilde{\psi}$ sur $X$,
on consid\ère $\tilde{\psi_{1}}$ et $\tilde{\psi_{2}}$
d\éfinies sur
$(\complex^{m+1}\backslash\bigcup_{i}\{z_{i}^{(0)}=0\})\times(\complex^{k+1}\backslash
\bigcup_{j}\{z_{j}^{(1)}=0\})$ par
$$\tilde{\psi_{1}}
=\ln\frac{(\mid z_{0}^{(0)}\mid...\mid
z_{k}^{(0)}\mid)^{\frac{2(m+1-k)}{k+1}}\times(\mid
z_{0}^{(1)}\mid...\mid z_{k}^{(1)}\mid)^{\frac{2k}{k+1}}}{(\mid
z_{0}^{(0)}\mid^{2}+...+\mid
z_{m}^{(0)}\mid^{2})^{(m+1-k)}\times(\mid
z_{0}^{(1)}\mid^{2}+...+\mid z_{k}^{(1)}\mid^{2})^{k}}
$$et $$\tilde{\psi_{2}}
=\ln\frac{(\mid z_{k+1}^{(0)}\mid...\mid
z_{m}^{(0)}\mid)^{\frac{2(m+1-k)}{m-k}}\times(\mid
z_{0}^{(1)}\mid...\mid z_{k}^{(1)}\mid)^{\frac{2k}{k+1}}}{(\mid
z_{0}^{(0)}\mid^{2}+...+\mid
z_{m}^{(0)}\mid^{2})^{(m+1-k)}\times(\mid
z_{0}^{(1)}\mid^{2}+...+\mid z_{k}^{(1)}\mid^{2})^{k}},$$ o\`{u}
$(z_{0}^{(0)},..,z_{m}^{(0)})$ sont les coordonn\ées de
$\complex^{m+1}$ et $(z_{0}^{(1)},..,z_{k}^{(1)})$ sont celles  de
$\complex^{k+1}$. Elles sont s\épar\ément homog\ènes de
degr\é z\éro en les composantes de chacun des vecteurs de
$\complex^{m+1}$ et $\complex^{k+1}$. Elles d\éfinissent alors
deux fonctions sur $\complex^{m+1}\times\complex^{k+1}$ et donc
sur $X$, par restriction. On pose $\tilde{\psi}=\inf
(\tilde{\psi_{1}},\tilde{\psi_{2}})$.

De m\^{e}me pour $Y$, on consid\ère les fonctions
\begin{eqnarray}
\hat{\psi_{1}}=\ln\{\frac{(\mid z_{0}^{(0)}\mid...\mid
z_{k}^{(0)}\mid)^{\frac{4}{k+1}}}{(\mid
z_{0}^{(0)}\mid^{2}+...+\mid z_{m}^{(0)}\mid^{2})^{2}} \times
\frac{(\mid z_{0}^{(1)}\mid...\mid
z_{k}^{(1)}\mid)^{\frac{2k}{k+1}}} {(\mid
z_{0}^{(1)}\mid^{2}+...+\mid z_{k}^{(1)}\mid^{2})^{k}}
\nonumber\\
\times\frac{(\mid z_{0}^{(2)}\mid...\mid
z_{m-k-1}^{(2)}\mid)^{\frac{2(m-k-1)}{m-k}}} {{(\mid
z_{0}^{(2)}\mid^{2}+...+\mid
z_{m-k-1}^{(2)}\mid^{2})^{(m-k-1)}}}\}\nonumber
\end{eqnarray}

et
\begin{eqnarray}
\hat{\psi_{2}}=\ln\{\frac{(\mid z_{k+1}^{(0)}\mid...\mid
z_{m}^{(0)}\mid)^{\frac{4}{m-k}}}{(\mid
z_{0}^{(0)}\mid^{2}+...+\mid z_{m}^{(0)}\mid^{2})^{2}} \times
\frac{(\mid z_{0}^{(1)}\mid...\mid
z_{k}^{(1)}\mid)^{\frac{2k}{k+1}}} {(\mid
z_{0}^{(1)}\mid^{2}+...+\mid z_{k}^{(1)}\mid^{2})^{k}}
\nonumber\\
\times\frac{(\mid z_{0}^{(2)}\mid...\mid
z_{m-k-1}^{(2)}\mid)^{\frac{2(m-k-1)}{m-k}}} {{(\mid
z_{0}^{(2)}\mid^{2}+...+\mid
z_{m-k-1}^{(2)}\mid^{2})^{(m-k-1)}}}\}.\nonumber
\end{eqnarray}
$\hat{\psi_{1}}$ et $\hat{\psi_{2}}$ sont deux fonctions sur
$$(\complex^{m+1}\backslash\bigcup_{i}\{z_{i}^{(0)}=0\})
\times(\complex^{k+1}\backslash\bigcup_{j}\{z_{j}^{(0)}=0\})
\times(\complex^{m-k}\backslash\bigcup_{q}\{z_{q}^{(0)}=0\})$$
o\`{u} $(z_{i}^{(0)})_{0\leq i\leq m}$, $(z_{j}^{(1)})_{0\leq
j\leq k}$ et $(z_{q}^{(2)})_{0\leq q\leq m-k-1}$ sont
respectivement les coordonn\ées sur $\complex^{m+1}$,
$\complex^{k+1}$ et $\complex^{m-k}$. Elles sont
s\épar\ément homog\ènes de degr\é z\éro en les
variables de $\complex^{m+1}$, $\complex^{k+1}$ et
$\complex^{m-k}$. Elles d\éfinissent des fonctions sur
$\proj_{m}\complex\times\proj_{k}\complex\times\proj_{m-k-1}\complex$,
et donc, par restriction, sur $Y$. On pose alors
$\hat{\psi}=\inf(\hat{\psi_{1}},\hat{\psi_{2}})$

Enon\c{c}ons \`a pr\'esent les
principaux r\ésultats de cet article.\\
\begin{theorem}\label{th1}
Soit $\varphi\in C^{\infty}(\proj_{m}\complex)$ une fonction
$g$-admissible et $G$-invariante, v\érifiant $\sup\varphi = 0$
sur $\proj_{m}\complex$. On a alors $\varphi \geq \psi$.
\end{theorem}
On en d\éduit le corollaire suivant.
\begin{coro}\label{coro1} Pour tout $\alpha <
\inf(\frac{k+1}{m+1},\frac{m-k}{m+1})$, on a l'in\égalit\é
de type H\H{o}rmander suivante (voir \cite{H} th. 4.4.5):
$$\int_{\proj_{m}\complex} \exp(-\alpha\varphi ) dv \leq Cst,$$
pour toute fonction $\varphi\in C^{\infty}(\proj_{m}\complex)$,
$g$-admissible, G-invariante, v\érifiant \\$\sup \varphi =0$
sur $\proj_{m}\complex$. $dv$ est l'\él\ément de volume sur
$\proj_{m}$ relatif \`{a} la m\étrique $g$.
\end{coro}
\begin{theorem}\label{th2}
Soit $\varphi\in C^{\infty}(X)$ une fonction
$\tilde{g}$-admissible et $G$-invariante, v\érifiant
$\sup\varphi = 0$ sur $X$. On a alors $\varphi \geq \tilde{\psi}$.
\end{theorem}
\begin{coro}\label{coro2} Pour tout $\alpha <
\inf(\frac{k+1}{m+1},\frac{m-k}{m-k+1})$, on a l'in\égalit\é
$$\int_{X} \exp(-\alpha\varphi ) d\tilde{v} \leq Cst,$$
pour toute fonction $\varphi\in C^{\infty}(X)$,
$\tilde{g}$-admissible, G-invariante, v\érifiant $\sup \varphi
=0$ sur $X$. $d\tilde{v}$ est l'\él\ément de volume sur $X$
relatif \`{a} la m\étrique $\tilde{g}$.
\end{coro}
\begin{theorem}\label{th3}
Soit $\varphi\in C^{\infty}(Y)$ une fonction $\hat{g}$-admissible
et $G$-invariante, v\érifiant $\sup\varphi = 0$ sur $Y$. On a
alors $\varphi \geq\hat{\psi}$.
\end{theorem}
\begin{coro}\label{coro3} Pour tout $\alpha <1/2$, on a
l'in\égalit\é
$$\int_{Y} \exp(-\alpha\varphi ) d\hat{v} \leq Cst,$$
pour toute fonction $\varphi\in C^{\infty}(Y)$,
$\hat{g}$-admissible, G-invariante, v\érifiant $\sup \varphi
=0$ sur $Y$. $d\hat{v}$ est l'\él\ément de volume sur $Y$
relatif \`{a} la m\étrique $\large\hat{g}$.
\end{coro}

{\section{Preuve des r\ésultats sur
$\proj_{m}\complex$.}\label{p}} {\subsection{Preuve du
th\éor\ème \ref{th1}.}} Pour le th\éor\ème \ref{th1},
on utilisera l'invariance des fonctions $\varphi
([z_{0},...,z_{m}])$ par le groupe $G$, afin de les
consid\érer, au lemme \ref{lem1}, comme des fonctions $\varphi
([1,x_{1},...,x_{m}])$ des variables r\'eelles $ x_{i} =\mid z_{i}
\mid$, $i\in\{1,..,m\}$, puis, au lemme \ref{lem2}, comme des
fonctions $\varphi ([x_{0},...,x_{k},1,x_{k+2},..,x_{m}])$ des
variables r\'eelles $ x_{i} =\mid z_{i} \mid$,
$i\in\{0,..,k,k+2,..,m\}$.
\begin{lem}\label{lem1}
Soit une fonction $\varphi\in C^{\infty}(\proj_{m}\complex)$,
$g$-admissible, G-invariante. \\Si $x_{i}=\mid z_{i} \mid >0$ pour
tout $i$, dans la carte $\{z_{0}\neq 0\}$, on a
\begin{eqnarray}\label{eq1}
(\varphi -\psi )([1,x_{1},..,x_{m}])\geq (\varphi -\psi
)([1,x_{1},..,x_{k};\zeta^{[m-k]}]),
\end{eqnarray}
o\`u $\zeta=(x_{k+1}...x_{m})^{1/(m-k)}$ et $\zeta^{[m-k]}=(\zeta
,..,\zeta)\in \complex^{m-k}$.
\end{lem}
{\bf Preuve.} La d\émonstration se fait par r\écurrence.
Supposons que pour $k+1\leq p < m$ et pour tout $(x_{1},..,x_{m})
\in\reel^{m}$ avec $x_{i}>0$ on ait
\begin{eqnarray}\label{eq2}
&&(\varphi -\psi )([1,x_{1},..,x_{m}])\geq \nonumber \\&& (\varphi
-\psi )([1,x_{1},..,x_{k};(x_{k+1}...x_{p})^{\frac{1}{p-k}},..,
(x_{k+1}...x_{p})^{\frac{1}{p-k}},x_{p+1},..x_{m}]).
\end{eqnarray}
Cette propri\ét\é est claire pour $p=k+1$. Si
l'in\égalit\é (\ref{eq2}) n'\était pas satisfaite au rang
$p+1$, il existerait alors un point
$(x_{1}^{0},..,x_{m}^{0})\in\reel^{m}$ avec $x_{i}^{0}>0$ pour
tout $i$, tel que
\begin{eqnarray}\label{eq3}
&&(\varphi -\psi )([1,x_{1}^{0},..,x_{m}^{0}])< \nonumber \\&&
(\varphi
-\psi)([1,x_{1}^{0},..,x_{k}^{0};(x_{k+1}^{0}...x_{p+1}^{0})^
{\frac{1}{p+1-k}},..,(x_{k+1}^{0}...x_{p+1}^{0})^{\frac{1}{p+1-k}},x_{p+2}^{0},..,x_{m}^{0}].
\end{eqnarray}

En utilisant la continuit\é de $(\varphi -\psi)$, on peut
supposer, quitte \`{a} en modifier l\ég\érement les
coordonn\ées, que le point $([1,x_{1}^{0},..,x_{m}^{0}])$ de
l'in\égalit\é (\ref{eq3}), v\érifie
$$(x_{1}^{0}..x_{k}^{0})^{1/(k+1)}\neq
(x_{k+1}^{0}..x_{m}^{0})^{1/(m-k)},$$propri\ét\é dont on
aura besoin plus loin. En utilisant la $G$-invariance de
$\varphi$, on peut supposer que $x_{k+1}^{0}\leq ...\leq
x_{m}^{0}$. D'autre part, en tenant encore compte de la $G$
invariance de $\varphi$ et de l'hypoth\èse de r\écurrence
(\ref{eq2}) en les points
$$[1,x_{1}^{0},..,x_{k}^{0};x_{k+1}^{0},x_{k+2}^{0},..,x_{p}^{0},
x_{p+1}^{0},x_{p+2}^{0},..,x_{m}^{0}]$$ et
$$[1,x_{1}^{0},..,x_{k}^{0};x_{k+2}^{0},
x_{k+3}^{0},..,x_{p}^{0},x_{p+1}^{0},x_{k+1}^{0},x_{p+2}^{0},..,x_{m}^{0}],$$
on peut \écrire
\begin{eqnarray}\label{eq4}&&(\varphi-\psi )([1,x_{1}^{0},..,x_{m}^{0}])\geq\nonumber\\
&& (\varphi-\psi
)([1,x_{1}^{0},..,x_{k}^{0};(x_{k+1}^{0}...x_{p}^{0})
^{\frac{1}{p-k}},..,(x_{k+1}^{0}...x_{p}^{0})^
{\frac{1}{p-k}},x_{p+1}^{0},x_{p+2}^{0},..,x_{m}^{0}]
\end{eqnarray}
et
\begin{eqnarray}\label{eq5}
&&(\varphi-\psi)([1,x_{1}^{0},..,x_{k}^{0};x_{k+2}^{0},..,x_{p+1}^{0},x_{k+1}^{0},
x_{p+2}^{0},..,x_{m}^{0}])\geq
\nonumber \\&&
(\varphi-\psi)([1,x_{1}^{0},..,x_{k}^{0};(x_{k+2}^{0}..x_{p+1}^{0})
^{\frac{1}{p-k}},..,(x_{k+2}^{0}..x_{p+1}^{0})^{\frac{1}{p-k}},x_{k+1}^{0},
x_{p+2}^{0},..,x_{m}^{0}].
\end{eqnarray}
Consid\érons maintenant la courbe $C$ d'\équation
$$t^{p-k}x=x_{k+1}^{0}...x_{p+1}^{0}$$ dans le plan r\éel
$\{[1,x_{1}^{0},..,x_{k}^{0},t,..,t,x,x_{p+2}^{0},..,x_{m}^{0}]\}$
param\étr\é par les variables $t$ et $x$. Les points
$$P_{1}=[1,x_{1}^{0},..,x_{k}^{0};(x_{k+1}^{0}...x_{p}^{0})^
{\frac{1}{p-k}},..,(x_{k+1}^{0}...x_{p}^{0})^{\frac{1}{p-k}},
x_{p+1}^{0},x_{p+2}^{0},..,x_{m}^{0}]$$ et
$$P_{2}=[1,x_{1}^{0},..,x_{k}^{0};(x_{k+2}^{0}...x_{p+1}^{0})
^{\frac{1}{p-k}},..,(x_{k+2}^{0}...x_{p+1}^{0})
^{\frac{1}{p-k}},x_{k+1}^{0},x_{p+2}^{0},..,x_{m}^{0}]$$
appartiennent \`{a} la courbe $C$.Notons que les r\éels
$x_{i}^{0}$ pour $k+1\leq i\leq p+1$ ne sont pas tous \égaux,
sinon (\ref{eq3})deviendrait une \égalit\é.\\ Par suite,
sachant que l'on a choisi $x_{k+1}^{0}\leq ...\leq x_{p+1}^{0}$,
les points distincts $P_{1}$ et $P_{2}$ se trouvent strictement de
part et d'autre de la diagonale $t=x$ du plan pr\éc\édent.
\\Or la courbe $C$ coupe cette diagonale en le point

$$P_{3}=[1,x_{1}^{0},..,x_{k}^{0};
(x_{k+1}^{0}...x_{p+1}^{0})^{\frac{1}{p+1-k}},..,
(x_{k+1}^{0}...x_{p+1}^{0})^{\frac{1}{p+1-k}},x_{p+2}^{0},..,x_{m}^{0}]$$
qui intervient dans l'in\égalit\é (\ref{eq3}). D'autre part,
en utilisant les relations (\ref{eq3}), (\ref{eq4}) et (\ref{eq5})
on obtient :$$ (\varphi - \psi)(P_{3})> (\varphi - \psi
)(P_{1})\mbox{ et }(\varphi - \psi )(P_{3})> (\varphi - \psi
)(P_{2}),$$ ce qui prouve que la fonction $(\varphi - \psi )$
admet un maximum local sur la courbe $C$. En cons\équence, la
restriction de la fonction $G$-invariante $(\varphi - \psi )$
\`{a} la courbe holomorphe (toujours not\'ee $C$) d'\équation
$\xi^{p-k}z=x_{k+1}^{0}...x_{p+1}^{0}$ du plan complexe
${[1,x_{1}^{0},..,x_{k}^{0};\xi,..,\xi,z,x_{p+2}^{0},..,x_{m}^{0}]}$
atteint un maximum local en un point $P=C(\zeta)$. Posons
$C(\zeta)=[1,C^{1}(\zeta),..,C^{m}(\zeta)]$ et
$\dot{C}^{\lambda}(\xi)=\frac{d C^{\lambda}}{d \xi}(\xi) \mbox{ et
} \dot{C}^{\overline{\mu}}(\xi)=\overline{\dot{C}^{\mu}(\xi)}.$

Sachant que l'on a choisi le point $[1,x_{1}^{0},..,x_{m}^{0}]$ de
sorte que $$(x_{1}^{0}...x_{k}^{0})^{1/(k+1)}\neq
(x_{k+1}^{0}...x_{m}^{0})^{1/(m-k)},$$ l'\équation de la
courbe $C$ et les d\éfinitions de $\psi_{1}$ et $\psi_{2}$ montrent 
qu'en tout point de $C$ 
$$ \psi_{1}([1,x_{1}^{0},..,x_{k}^{0};\xi,..,\xi,z,
x_{p+2}^{0},..,x_{m}^{0}])\neq
\psi_{2}([1,x_{1}^{0},..,x_{k}^{0};\xi,..,\xi,z,x_{p+2}^{0},..,x_{m}^{0}]).$$
On peut alors supposer que $\psi=\psi_{1}$ dans un voisinage de $P$, 
la preuve \étant identique si l'on suppose $\psi=\psi_{2}$ dans
ce voisinage. On a donc :

$$\frac {\partial^{2}}{\partial \xi\partial \overline{\xi}}
\{(\varphi -\psi_{1})(C(\zeta))\} =\frac{\partial^{2}(\varphi
-\psi_{1})}{\partial z_{\lambda}
\partial \overline{z}_{\mu} }(C(\zeta))\dot{C}^{\lambda}(\zeta)
\dot{C}^{\overline{\mu}}(\zeta)$$ est n\égatif ou nul.
Comme $-\frac{\partial^{2}\psi_{1}}{\partial z_{\lambda}
\partial \overline{z}_{\mu}}= g_{\lambda\overline{\mu}},$
ceci exprime que la forme hermitienne de matrice:
$$(g_{\lambda\overline{\mu}}+\frac{\partial^{2}\varphi}{\partial
z_{\lambda}
\partial \overline{z}_{\mu}})_{\lambda ,\mu}=(\frac{\partial^{2}(\varphi -\psi_{1} )}
{\partial z_{\lambda}
\partial \overline{z}_{\mu}})_{\lambda,\mu}$$
est n\égative en $P = C(\zeta)$. On en
d\éduit une contradiction avec la $g$-admissibilit\é de
$\varphi$ en $P$. D'o\`u l'in\'egalit\'e (\ref{eq2}) au rang $p+1$ et, par
cons\'equent, le lemme \ref{lem1}.

\begin{lem}\label{lem2}
Soit une fonction $\varphi\in
C^{\infty}(\proj_{m}\complex)$, $g$-admissible, G-invariante.// Si
$x_{i}=\mid z_{i}\mid >0$ pour tout $i$, dans la carte
$\{z_{k+1}\neq 0\}$ on a
\begin{eqnarray}\label{eq6}
(\varphi -\psi )([x_{0},x_{1},..,x_{k};1,x_{k+2},..,x_{m}])\geq
(\varphi -\psi )([\eta,\eta,..,\eta;1,x_{k+2},..,x_{m}]),
\end{eqnarray}
o\`u $\eta=(x_{0}x_{1}...x_{k})^{1/(k+1)}$.
\end{lem}
{\bf Preuve.} Comme dans le lemme \ref{lem1}, la preuve s'effectue
par r\écurrence. Supposons que pour $0\leq p < k$ et
pour tout $(x_{0},..,x_{k};x_{k+2},..,x_{m}) \in\reel^{m}$ avec
$x_{i}>0$, on ait
\begin{eqnarray}\label{eq7}
&&(\varphi -\psi )([x_{0},..,x_{k},1,x_{k+2},..,x_{m}])\geq
\nonumber
\\&& (\varphi -\psi
)([(x_{0}...x_{p})^{\frac{1}{p+1}},..,
(x_{0}...x_{p})^{\frac{1}{p+1}},x_{p+1},..,x_{k};1,x_{k+2},..,x_{m}]).
\end{eqnarray}
Cette hypoth\èse est v\érifi\ée pour $p=0$. Si
l'in\égalit\é (\ref{eq7}) n'\était pas satisfaite au rang
$p+1$, il existerait un point
$(x_{0}^{0},..,x_{k};x_{k+2},..,x_{m}^{0})\in\reel^{m}$ avec
$x_{i}^{0}>0$ pour tout $i$, tel que :
\begin{eqnarray}\label{eq8}
&&(\varphi -\psi
)([x_{0}^{0},..,x_{k}^{0};1,x_{k+2}^{0},..,x_{m}^{0}])< \nonumber
\\&& (\varphi
-\psi)([(x_{0}^{0}..x_{p+1}^{0})^{\frac{1}{p+2}},..,(x_{0}^{0}..x_{p+1}^{0})^{\frac{1}{p+2}},
x_{p+2}^{0},...,x_{k}^{0};1,x_{k+2}^{0},..,x_{m}^{0}]).
\end{eqnarray}
Comme au lemme \ref{lem1}, on peut supposer que le point
$[x_{0}^{0},..,x_{k}^{0};1,x_{k+2}^{0},..,x_{m}^{0}]$ v\érifie
$$(x_{0}^{0}..x_{k}^{0})^{1/(k+1)}\neq
(x_{k+2}^{0}..x_{m}^{0})^{1/(m-k)},$$ et que $x_{0}^{0}\leq
...\leq x_{p+1}^{0}$. D'autre part, en tenant compte de la $G$-invariance 
de $\varphi$ et de l'hypoth\èse de r\écurrence
(\ref{eq7}) en les points 
$$[x_{0}^{0},x_{1}^{0},..,x_{p}^{0},x_{p+1}^{0},..,x_{k}^{0};1,x_{k+2}^{0},..,x_{m}^{0}]$$
et
$$[x_{1}^{0},..,x_{p+1}^{0},x_{0}^{0},..,x_{k}^{0};1,x_{k+2}^{0},..,x_{m}^{0}],$$
on a 
\begin{eqnarray}\label{eq9}
&&(\varphi-\psi )([x_{0}^{0},..,x_{k}^{0};1,x_{k+2}^{0},
..,x_{m}^{0}])\geq\nonumber\\&& (\varphi-\psi
)([(x_{0}^{0}..x_{p}^{0})^{\frac{1}{p+1}},..,(x_{0}^{0}..x_{p}^{0})^{\frac{1}{p+1}},
x_{p+1}^{0},..,x_{k}^{0};1,x_{k+2}^{0},..,x_{m}^{0}]) 
\end{eqnarray}
et
\begin{eqnarray}\label{eq10}
&&(\varphi-\psi)([x_{1}^{0},..,x_{p+1}^{0},x_{0}^{0},..,x_{k}^{0};1,x_{k+2}^{0},..,x_{m}^{0}])\geq
\nonumber \\&&
 (\varphi-\psi)([(x_{1}^{0}..x_{p+1}^{0})^{\frac{1}{p+1}},..,(x_{1}^{0}..x_{p+1}^{0})^
{\frac{1}{p+1}},
x_{0}^{0},x_{p+2}^{0},..,x_{k}^{0};1,x_{k+2}^{0},..,x_{m}^{0}]).
\end{eqnarray}
Consid\érons maintenant la courbe $C$ d'\équation
$$t^{p+1}x=x_{0}^{0}...x_{p+1}^{0}$$ du plan r\éel
$\{[t,..,t,x,x_{p+2}^{0},..,x_{k}^{0};1,x_{k+2}^{0},..,x_{m}^{0}]\}$
param\étr\é par les variables $t$ et $x$. Les points 
$$Q_{1}=[(x_{0}^{0}..x_{p}^{0})^{\frac{1}{p+1}},..,
(x_{0}^{0}..x_{p}^{0})^{\frac{1}{p+1}},x_{p+1}^{0},x_{p+2}^{0},..,x_{k}^{0};1,
x_{k+2}^{0},..,x_{m}^{0}]$$ et
$$Q_{2}=[(x_{1}^{0}..x_{p+1}^{0})^{\frac{1}{p+1}},..,(x_{1}^{0}..x_{p+1}^{0})^{\frac{1}{p+1}},
x_{0}^{0},x_{p+2}^{0},..,x_{k}^{0};1,x_{k+2}^{0},..,x_{m}^{0}]$$
appartiennent \`{a} la courbe $C$.\\D'autre part 
 les r\éels $x_{i}^{0}$ pour $0\leq i\leq p+1$ ne sont pas tous
\égaux, sinon (\ref{eq8}) serait une \égalit\é.\\
Par suite, sachant que l'on a choisi $x_{0}^{0}\leq ...\leq
x_{p+1}^{0}$, les points distincts $Q_{1}$ et $Q_{2}$ se trouvent strictement 
de part et
d'autre de la diagonale $t=x$ du plan pr\éc\édent. \\Or la
courbe $C$ coupe cette diagonale en le point 

$$Q_{3}=[(x_{0}^{0}...x_{p+1}^{0})^{\frac{1}{p+2}},..,
(x_{0}^{0}...x_{p+1}^{0})^{\frac{1}{p+2}},x_{p+2}^{0},..,x_{k}^{0};
1,x_{k+2}^{0},..,x_{m}^{0}]$$ qui intervient dans
l'in\égalit\é (\ref{eq8}). D'autre part, les relation
(\ref{eq8}), (\ref{eq9}) et (\ref{eq10}) donnent 
$$ (\varphi -
\psi)(Q_{3})> (\varphi - \psi )(Q_{1})\mbox{ et }(\varphi - \psi
)(Q_{3})> (\varphi - \psi )(Q_{2}),$$ ce qui prouve que la
fonction $(\varphi - \psi )$ admet un maximum local sur la courbe
$C$.\\
Sachant que l'on a choisi le point
$[x_{0}^{0},..,x_{k}^{0};1,x_{k+2}^{0},..,x_{m}^{0}]$ de sorte que
$$(x_{0}^{0}..x_{k}^{0})^{1/(k+1)}\neq
(x_{k+2}^{0}..x_{m}^{0})^{1/(m-k)},$$
on conclut de la
m\^{e}me mani\ère qu'au lemme pr\éc\édent en consid\'{é}rant la restriction de 
$(\varphi -\psi)$ \`{a} une courbe holomorphe convenable.
\begin{lem}\label{coro4}
Etant donn\ée une fonction $\varphi\in
C^{\infty}(\proj_{m}\complex)$, $g$-admissible, G-invariante, si 
$x_{i}=\mid z_{i}\mid >0$ pour tout $i$, on a 
\begin{eqnarray}\label{eq11}
(\varphi -\psi )([1,x_{1},..,x_{k};x_{k+1},x_{k+2},..,x_{m}])\geq
(\varphi -\psi )([1^{[k+1]};\nu^{[m-k]}]),
\end{eqnarray}
o\`u
$\nu=(x_{1}...x_{k})^{-1/(k+1)}(x_{k+1}...x_{m})^{1/(m-k)}$.
\end{lem}
{\bf Preuve}. D'apr\ès le lemme \ref{lem1}, on a 
\begin{eqnarray*}
   & & (\varphi -\psi
)([1,x_{1},..,x_{k};x_{k+1},x_{k+2},..,x_{m}])\nonumber\\ \\
  &\geq& (\varphi -\psi
)([1,x_{1},..,x_{k};(x_{k+1}...x_{m})^{1/(m-k)},..,(x_{k+1}...x_{m})^{1/(m-k)}])\nonumber \\
  &=& (\varphi -\psi )([\frac{1}{(x_{k+1}...x_{m})^{1/(m-k)}},
\frac{x_{1}}{(x_{k+1}...x_{m})^{1/(m-k)}},..,
\frac{x_{k}}{(x_{k+1}...x_{m})^{1/(m-k)}};1,..,1]).\nonumber
\end{eqnarray*}
La $(k+2)$-i\ème composante homog\ène de ce dernier point \étant
\égale \`{a} $1$, le 
lemme \ref{lem2} nous permet d'\écrire 
\begin{eqnarray*}
   && (\varphi -\psi
)([1,x_{1},..,x_{k};x_{k+1},x_{k+2},..,x_{m}])\nonumber \\
   &\geq& (\varphi -\psi )([\frac{1}{(x_{k+1}...x_{m})^{1/(m-k)}},
\frac{x_{1}}{(x_{k+1}...x_{m})^{1/(m-k)}},..,
\frac{x_{k}}{(x_{k+1}...x_{m})^{1/(m-k)}};1,..,1])\nonumber \\
   &\geq& (\varphi -\psi
)([\frac{(x_{1}...x_{k})^{1/(k+1)}}{(x_{k+1}...x_{m})^{1/(m-k)}},..,
\frac{(x_{1}...x_{k})^{1/(k+1)}}{(x_{k+1}...x_{m})^{1/(m-k)}};1,..,1])\nonumber \\
  &=& (\varphi -\psi
)([1,..,1;\frac{(x_{k+1}...x_{m})^{1/(m-k)}}{(x_{1}...x_{k})^{1/(k+1)}},..,
\frac{(x_{k+1}...x_{m})^{1/(m-k)}}{(x_{1}...x_{k})^{1/(k+1)}}]),\nonumber
\end{eqnarray*}
d'o\`{u} la minoration (\ref{eq11}).

Pour la suite de la preuve du th\éor\ème \ref{th1},
pr\écisons que la $G$ invariance ne nous permet pas d'aller
plus loin, contrairement \`{a} ce qui a \ét\é fait dans
\cite{B1} o\`{u} le groupe d'automorphismes est plus gros. 
C'est pour cette raison que l'on ne peut
pas conserver ici la fonction extr\émale de \cite{B1}. Celle qui apparaît 
dans cet article nous permet de passer directement du
lemme \ref{coro4} \`{a} la derni\ère \étape, \`{a}
savoir le lemme suivant :
\begin{lem}\label{lem4}
Etant donn\ée une fonction $\varphi\in
C^{\infty}(\proj_{m}\complex)$, $g$-admissible, G-invariante 
et telle que $\sup=0$ sur $\proj_{m}\complex$, alors pour tout $\zeta >0$, on a 
\begin{eqnarray}\label{eq12}
(\varphi -\psi )([1^{[k+1]};\zeta^{[m-k]}])\geq 0,
\end{eqnarray}
\end{lem}
{\bf Preuve.} On raisonne sur la position du point $P_{0}\in\proj_{m}\complex$
o\`{u} $\varphi$ atteint son maximum. En vertu de la $G$-invariance
de $\varphi$, on peut supposer qu'il s'\écrit sous la forme 
$$P_{0}=[y^{0}_{0},..,y^{0}_{k};y^{0}_{k+1},..,y^{0}_{m}],$$
o\`{u} les $y^{0}_{i}$ sont des r\éels positifs v\érifiant 
$y^{0}_{0}\geq y^{0}_{1}\geq..\geq y^{0}_{k}$ et $y^{0}_{k+1}\geq
y^{0}_{k+2}\geq..\geq y^{0}_{m}$. Deux cas se pr\ésentent : ou
bien l'un des $y^{0}_{0},..,y^{0}_{k}$ est non nul, ou bien tous
les $y^{0}_{0},..,y^{0}_{k}$ sont nuls.
\\{\bf{Cas A}} :
l'un des $y^{0}_{0},..,y^{0}_{k}$ est non nul. On peut alors se
placer dans la carte $\{z_{0}\neq 0\}$ et \écrire le point
$P_{0}$ sous la forme 
$$P_{0}=[1,x^{0}_{1},..,x^{0}_{k};x^{0}_{k+1},..,x^{0}_{m}],$$
o\`{u} les r\éels positifs $x^{0}_{i}$ v\érifient :
$1\geq x^{0}_{1}\geq..\geq x^{0}_{k}$ et $x^{0}_{k+1}\geq ...\geq
x^{0}_{m}$.

Raisonnons par l'absurde, et supposons qu'il existe un point
$P_{1}=[1^{[k+1]};\zeta_{0}^{[m-k]}]$ tel que
l'on ait $\zeta_{0}>0$ et 
\begin{equation}\label{p1}
(\varphi -\psi )(P_{1})<0.
\end{equation}
 On envisage alors les deux sous-cas suivants : $x^{0}_{k+1}<\zeta_{0}$ 
puis $x^{0}_{k+1}\geq \zeta_{0}$.
\begin{itemize}
    \item \underline{$x^{0}_{k+1}<\zeta_{0}$}.
\end{itemize}
On introduit alors la fonction auxiliaire
$$\psi_{0}=\log\frac{\mid z_{0}\mid^{2(m+1)}}{
    (\mid z_{0}\mid^{2}+...+\mid z_{m}\mid^{2})^{m+1}}.$$
D'une part, puisque $\varphi\leq 0$,
\begin{equation}\label{p}
(\varphi -\psi_{0})([1,0,...,0])=\varphi ([1,0,...,0])\leq 0.
\end{equation}
De plus, sachant que $\varphi (P_{0})=0$ et $\psi_{0}\leq 0$, on a
\begin{equation}\label{p3}
(\varphi -\psi_{0})(P_{0})\geq 0.
\end{equation}
Si $P_{0}\neq [1,0,..,0]$, $\psi_{0}(P_{0})< 0$ et l'in\égalit\é (\ref{p3}) est
alors stricte. Si $P_{0}= [1,0,..,0]$, 
quitte \`{a} se placer en un point $P$ arbitrairement voisin de
$P_{0}$, on peut supposer $(\varphi -\psi_{0})(P)> 0$. En effet, si dans
un voisinage de $P_{0}$ on avait $(\varphi -\psi_{0})\leq 0$, comme 
$(\varphi -\psi_{0})(P_{0})= 0$, $(\varphi -\psi_{0})$ admettrait
alors un maximun local en $P_{0}$, ce qui mettrait en
d\éfaut l'admissibilit\é de $\varphi$ en ce point, sachant que 
$$\partial_{\lambda\overline{\mu}}(\varphi -\psi_{0})(P_{0})=
(g_{\lambda\overline{\mu}}+\partial_{\lambda\overline{\mu}}\varphi )(P_{0}).$$ Dans tous les cas, 
on peut
donc affirmer qu'il existe un point $P'_{0}=[1,{a}_{1},..,{a}_{m}]$
v\érifiant
\begin{equation}\label{p3'}
(\varphi -\psi_{0})(P'_{0})> 0.
\end{equation}
Par continuit\é et $G$-invariance de $\varphi$, on peut supposer $1>a_{1}\geq...\geq a_{k}>0$ et 
$\zeta_{0}> a_{k+1}\geq...\geq a_{m}>0$.
D'autre part, l'in\égalit\é (\ref{p1}) jointe aux d\'efinitions de $P_1$, $\psi_0$, $\psi_1$ 
et $\psi=\inf (\psi_1 , \psi_2)$ implique
\begin{equation}\label{p1'}
(\varphi -\psi_{0})(P_{1})=(\varphi -\psi_{1})(P_{1})\leq(\varphi
-\psi)(P_{1})< 0.
\end{equation}

La courbe :
$$[0,1]\ni t\rightarrow [1,t,t^{(\ln a_{2})/(\ln a_{1})},..,t^{(\ln a_{k})/(\ln a_{1})};
\zeta_{0}t^{\frac{\ln(a_{k+1}/\zeta_{0})}{\ln
a_{1}}},..,\zeta_{0}t^{\frac{\ln(a_{m}/\zeta_{0})}{\ln
a_{1}}}]$$ passe par $[1,0,..,0]$ en $t=0$ puis par ${P'}_{0}$
en $t=a_{1}$ et enfin par le point $P_{1}$ en $t=1$, valeurs
en lesquelles, d'apr\ès (\ref{p}), (\ref{p3'}) et (\ref{p1'}),
$(\varphi-\psi_{0})$ est respectivement n\égative, positive
puis \`{a} nouveau n\égative. L'invariance de cette fonction
par l'action des $\exp(i\theta)$, permet donc de d\éduire que
$(\varphi-\psi_{0})$ atteint un maximum sur la courbe holomorphe,
complexifi\ée de la courbe d\écrite plus haut, ce qui
contredit encore une fois l'admissibilit\é de $\varphi$.
\begin{itemize}
\item \underline{$x^{0}_{k+1}\geq\zeta_{0}$}.
\end{itemize}
D\ésignons dans ce cas par $p\in \{1,..,m-k\}$ l'entier pour
lequel on a $$x^{0}_{k+1}\geq ...\geq x^{0}_{k+p}>\zeta_{0} \mbox{
et } \zeta_{0}\geq x^{0}_{k+p+1}\geq ...\geq x^{0}_{m},$$ et
consid\érons la fonction auxiliaire

$$\psi_{k+1}=\log\frac{\mid z_{k+1}\mid^{2(m+1)}}{
    (\mid z_{0}\mid^{2}+...+\mid z_{m}\mid^{2})^{m+1}}.$$
On a
\begin{equation}\label{p3''}
(\varphi -\psi_{k+1})(P_{0})> 0.
\end{equation}
La fonction $(\varphi -\psi_{k+1})$ \étant continue, quitte
\`{a} se placer en un point voisin de $P_{0}$, on peut supposer
tous les $x_{i}^{0}$ non nuls. Posons alors:
$$\alpha_{2}=\frac{\ln x_{2}^{0}}{\ln x_{1}^{0}},..,\alpha_{k}=\frac{\ln x_{k}^{0}}{\ln x_{1}^{0}};
\alpha_{k+1}=\frac{\ln(x_{k+1}^{0}/\zeta_{0})}{\ln x_{1}^{0}},..,
\alpha_{m}=\frac{\ln(x_{m}^{0}/\zeta_{0})}{\ln x_{1}^{0}}.$$
Sachant que l'on a $1\geq x_{1}^{0}\geq..\geq x_{k}^{0}$;
$x_{k+1}^{0}\geq..\geq x_{k+p}^{0}\geq\zeta_{0}$ et $\zeta_{0}\geq
x_{k+p+1}^{0}\geq..\geq x_{m}^{0}$, on en d\éduit que 
$\alpha_{2},..,\alpha_{k}\geq 0$, $\alpha_{k+1}\leq ...\leq
\alpha_{p+k}\leq 0$ et $\alpha_{p+k+1},..,\alpha_{m}\geq 0$, d'o\`{u}, en notant 
$P_{\varepsilon}=[1,\varepsilon,\varepsilon^{\alpha_{2}},..,\varepsilon^{\alpha_{k}},
\zeta_{0}\varepsilon^{\alpha_{k+1}},..,\zeta_{0}\varepsilon^{\alpha_{m}}]$, on a 
\begin{eqnarray*}
\lim_{\varepsilon\rightarrow
0}\psi_{k+1}(P_{\varepsilon }) 
&=& \lim_{\varepsilon\rightarrow
0}\ln\frac{\zeta_{0}^{2(m+1)}\varepsilon^{2\alpha_{k+1}(m+1)}}
{[1+\varepsilon^{2}+\varepsilon^{2\alpha_{2}}+..+\varepsilon^{2\alpha_{k}}+
\zeta_{0}^{2}(\varepsilon^{2\alpha_{k+1}}+..+\varepsilon^{2\alpha_{m}})]^{m+1}}  \\
&=&\ln\lim_{t\rightarrow
\infty}\frac{\zeta_{0}^{2(m+1)}t^{-2\alpha_{k+1}(m+1)}}
{[\zeta_{0}^{2}(t^{-2\alpha_{k+1}}+t^{-2\alpha_{k+2}}+..+t^{-2\alpha_{p}})]^{(m+1)}}
= \ln 1=0,
\end{eqnarray*} 
$-\alpha_{k+1}$ \étant la plus grande puissance intervenant au
d\énominateur. Sachant que $\varphi ([P_{\varepsilon}])\leq 0$ et compte tenu de (\ref{p3''}) il
existe $\varepsilon_{0}$ tel que l'on ait
\begin{equation}\label{p4}
(\varphi
-\psi_{k+1})(P_{\varepsilon_{0} })\leq -\psi_{k+1}(P_{\varepsilon_{0} })<
(\varphi -\psi_{k+1})(P_{0}).
\end{equation}
D'autre part, l'in\'egalit\'e (\ref{p1}), jointe aux d\'efinitions de $P_1$, $\psi_{k+1}$, $\psi_{2}$ 
et $\psi = \inf (\psi_{1},\psi_{2})$ donne~: 
\begin{equation}\label{p5}
(\varphi
-\psi_{k+1})(P_{1})= (\varphi-\psi_{2})(P_{1})\leq
(\varphi -\psi)(P_{1})<0.
\end{equation}
La courbe
$$[\varepsilon_{0},1]\ni t\rightarrow [1,t,t^{\alpha_{2}},..,t^{\alpha_{k}},
\zeta_{0}t^{\alpha_{k+1}},..,\zeta_{0}t^{\alpha_{m}}],$$ passe par
$P_{\varepsilon_{0}}$
en $t=\varepsilon_{0}$ puis par $P_{0}$ en $t=x_{1}^{0}$ et enfin
par $P_{1}$ en $t=1$, ce qui, en vertu de 
(\ref{p4}), (\ref{p3''}) et (\ref{p5}) prouve l'existence d'un maximum local pour la 
fonction $(\varphi-\psi_{k+1})$ sur la
courbe pr\écit\ée. Ceci contredit, \`{a} l'instar du cas
pr\éc\édent, l'hypoth\èse d'admissibilit\é de la
fonction $\varphi$.

{\bf{Cas B :}}  $y^{0}_{0}=...=y^{0}_{k}=0$. On se
place alors dans la carte $\{z_{k+1}\neq 0\}$, de sorte que
le point $P_{0}$ o\`u $\varphi$ atteint son maximum \'egal \`a z\'ero puisse 
s'\'ecrire sous la forme 
$$P_{0}=[0,0,..,0;1,x_{k+2}^{0},..,x_{m}^{0}].$$
On peut aussi supposer, en utilisant la $G$-invariance de $\varphi$, que 
$1\geq x_{k+2}^{0}\geq..\geq x_{m}^{0}$. On montrera une version \'equivalente du lemme \ref{lem4}, \`a savoir que
\begin{equation}\label{equ0}
(\varphi -\psi )([\zeta^{[k+1]};1^{[m-k]}])\geq 0
\end{equation}
pour tout $\zeta >0$. 
Raisonnons par l'absurde, et supposons qu'il existe un point
$P_{k+1}=[\zeta_{0}^{[k+1]};1^{[m-k]}]$ tel que
l'on ait $\zeta_{0}>0$ et 
\begin{equation}\label{equ1}
(\varphi -\psi )(P_{k+1})<0.
\end{equation}
On consid\`ere alors la fonction auxiliaire
$\psi_{k+1}$ introduite plus haut.
Sachant que $\varphi\leq 0$, on a 
\begin{equation}\label{equ2}
(\varphi -\psi_{k+1})([0^{[k+1]};1,0,..,0])=\varphi ([0^{[k+1]};1,0,..,0])\leq 0.
\end{equation}
D'autre part, sachant que $\varphi (P_{0})=0$ et $\psi_{k+1}\leq 0$, on a
\begin{equation}
(\varphi -\psi_{k+1})(P_{0})=-\psi_{k+1}(P_{0})\geq 0.
\end{equation}
Cette in\'egalit\'e est stricte d\`es que $P_{0}\neq [0^{[k+1]};1,0,..,0]$. Si 
$P_{0} = [0^{[k+1]};1,0,..,0]$, quitte \`a se placer en un point arbitrairement voisin de 
$P_{0}$, on peut supposer la derni\`ere in\'egalit\'e stricte. En effet, si dans un voisinage de $P_{0}$, on avait 
$\varphi -\psi_{k+1}\leq 0$, alors $\varphi -\psi_{k+1}$ admettrait un maximum local en  
$P_{0}$, ce qui contredirait l'admissibilit\'e de $\varphi$ en $P_{0}$. A l'instar du cas A, 
il existe donc un point  $P'_{0}=[{c}_{0},..,{c}_{k};1,{c}_{k+2},..,{c}_{m}]$
v\érifiant
\begin{equation}\label{equ3}
(\varphi -\psi_{k+1})(P'_{0})> 0.
\end{equation}
Par continuit\é et $G$-invariance de $\varphi$, on peut supposer 
$\zeta_{0}> c_{0}\geq...\geq c_{k}>0$ et $1>c_{k+2}\geq...\geq c_{m}>0$. 
D'autre part, l'in\égalit\é (\ref{equ1}) jointe aux d\'efinitions de $P_{k+1}$, $\psi_{k+1}$, $\psi_2$ 
et $\psi=\inf (\psi_1 , \psi_2)$ implique
\begin{equation}\label{equ4}
(\varphi -\psi_{k+1})(P_{k+1})=(\varphi -\psi_{2})(P_{k+1})\leq(\varphi
-\psi)(P_{k+1})< 0.
\end{equation}

La courbe :
$$[0,1]\ni t\rightarrow [\zeta_{0}t^{\frac{\ln(c_{0}/\zeta_{0})}{\ln
c_{k+2}}},..,\zeta_{0}t^{\frac{\ln(c_{k}/\zeta_{0})}{\ln
c_{k+2}}};1,t,t^{(\ln c_{k+3})/(\ln c_{k+2})},..,t^{(\ln c_{m})/(\ln c_{k+2})}]$$ 
passe par $[0^{[k+1]},1,0,..,0]$ en $t=0$ puis par ${P'}_{0}$
en $t=c_{k+2}$ et enfin par le point $P_{k+1}$ en $t=1$, valeurs
en lesquelles, d'apr\ès (\ref{equ2}), (\ref{equ3}) et (\ref{equ4}),
$(\varphi-\psi_{k+1})$ est respectivement n\égative, positive
puis n\égative. L'invariance de cette fonction
par l'action des $\exp(i\theta)$, permet donc de d\éduire que
$(\varphi-\psi_{0})$ atteint un maximum sur la courbe holomorphe,
d\'eduite de la courbe d\écrite plus haut, ce qui
contredit l'admissibilit\é de $\varphi$, d'o\`u (\ref{equ0}) 
et le lemme \ref{lem4}.

{\subsection{Preuve du corollaire \ref{coro1}.}}

Soit $\varphi\in C^{\infty}(\proj_{m}\complex )$ une fonction
$g$-admissible et $G$-invariante, dont le $\sup$ sur
$\proj_{m}\complex$ est nul. D'apr\`es le th\'eor\`eme \ref{th1},
on a $\varphi \geq \psi$ et par suite, pour tout $\alpha \geq 0$,
$$\int_{\proj_{m}\complex} \exp(-\alpha\varphi ) dv \leq
\int_{\proj_{m}\complex}\exp(-\alpha\psi )dv .$$ Cherchons les
valeurs de $\alpha$ pour lesquelles cette derni\`ere int\'egrale
converge. Pour ce faire, on estimera
$\int_{\proj_{m}\complex}\exp(-\alpha\psi_{1} )dv$ et
$\int_{\proj_{m}\complex}\exp(-\alpha\psi_{2} )dv$ dans la carte
dense d\'efinie par $\{z_{0}=1\}$. Dans cette carte, l'\'el\'ement
de volume est donn\'e par 
$$dv=(i)^{m}\frac{dz_{1}\wedge d\overline{z}_{1}\wedge ...\wedge
dz_{m}\wedge d\overline{z}_{m}} {(1+\mid z_{1}\mid^{2}+...+\mid
z_{m}\mid^{2})^{m+1}}.$$ En utilisant la d\'efinition de $\psi_{1}$ qui
ne d\'epand que des $\mid z_{p}\mid$, par le changement
de variables $u_{p}=\mid z_{p}\mid^{2}$ on obtient~:
\begin{eqnarray}
\int_{\proj_{m}\complex}\exp(-\alpha\psi_{1} )dv  = Cst
\int_{0}^{+\infty}...\int_{0}^{+\infty}\frac{
(u_{1}...u_{k})^{-\alpha(m+1)/(k+1)}du_{1}...du_{m}} {(1+
u_{1}+...+u_{m})^{(1-\alpha)(m+1)}}.\nonumber
\end{eqnarray}
Cette int\égrale converge en z\éro si et seulement si~: 
$1-\frac{\alpha(m+1)}{k+1}>0$, 
c'est \`a dire $\alpha < (k+1)/(m+1)$. 
En l'infini, le passage en
coordonn\ées polaires ram\ène l'\étude \`{a} la
convergence de
$$\int_{a>0}^{\infty}r^{\frac{-\alpha k(m+1)}{k+1}}r^{(\alpha-1)(m+1)}
r^{m-1}dr,$$ ce qui donne la condition : $$\frac{-\alpha k(m+1)}{k+1}+
(\alpha-1)(m+1)+(m-1)+1<0,$$ et donc encore $\alpha < (k+1)/(m+1)$.
Ainsi 
$\int_{\proj_{m}\complex}\exp(-\alpha\psi_{1} )dv$ converge si et seulement si $\alpha <
\frac{k+1}{m+1}.$ D'autre part,
\begin{eqnarray}
\int_{\proj_{m}\complex}\exp(-\alpha\psi_{2} )dv  = Cst
\int_{0}^{+\infty}..\int_{0}^{+\infty}\frac{
(u_{k+1}..u_{m})^{-\alpha(m+1)/(m-k)}du_{1}..du_{m}}
{(1+u_{1}+..+u_{k}+u_{k+1}+..+u_{m})^{(1-\alpha)(m+1)}}.\nonumber
\end{eqnarray}
La convergence en z\éro exige $\frac{\alpha(m+1)}{m-k}<1$,
c'est \`{a} dire $\alpha < \frac{m-k}{m+1}$. En l'infini, le
passage en coordonn\ées polaires ram\ène l'\étude \`{a}
la convergence de $$\int_{a>0}^{\infty}r^{-\alpha
(m+1)}r^{(\alpha-1)(m+1)} r^{m-1}dr=\int_{a>0}^{\infty}r^{-2}dr,$$ int\égrale convergente
quelle que soit la valeur de $\alpha$. D'o\`{u} la
convergence de $\int_{\proj_{m}\complex}\exp(-\alpha\psi_{2} )dv$
pour $\alpha < \frac{m-k}{m+1}$ et le corollaire \ref{coro1}.

{\subsection{Preuve du th\éor\ème \ref{th2}.}} Comme pour le
th\éor\ème \ref{th1}, on utilisera l'invariance par le
groupe $G$ d\éfini dans l'introduction, des fonctions $\varphi
([z_{0},..,z_{k},z_{k+1},..,z_{m}],[\zeta_{0},..,\zeta_{k}])$
o\`{u} $(z_{0},..,z_{k})$ et $(\zeta_{0},..,\zeta_{k})$ sont
colin\éaires. Cel\`{a} nous permettra de les consid\érer,
dans le lemme \ref{lem21} comme des fonctions $\varphi
([1,x_{1},...,x_{m}],[1,x_{1},...,x_{k}])$ des variables r\'eelles
$ x_{i} =\mid z_{i} \mid >0$, $i\in\{1,..,m\}$, puis, dans le lemme
\ref{lem22} comme des fonctions $$\varphi
([x_{0},...,x_{k},1,x_{k+2},..,x_{m}],[x_{0},...,x_{k}])$$ des
variables r\'eelles $ x_{i} =\mid z_{i} \mid >0 $, $i\in\{0,..,k,k+2,..,m\}$.
Notons que ces rep\érages ne contiennent pas l'\éclatement, les $x_{i}$ \étant non nuls 
dans l'\énonc\é de ces deux lemmes.
\begin{lem}\label{lem21}
Etant donn\ée une fonction $\varphi\in C^{\infty}(X)$,
$\tilde{g}$-admissible, G-invariante. Si $x_{i}=\mid z_{i}\mid >0$ pour
tout $i$:
\begin{eqnarray}\label{eq21}
(\varphi -\tilde{\psi}
)([1,x_{1},..,x_{m}],[1,x_{1},..,x_{k}])\geq (\varphi
-\tilde{\psi}
)([1,x_{1},..,x_{k};\zeta^{[m-k]}],[1,x_{1},..,x_{m}]),
\end{eqnarray}
o\`u $\zeta^{[m-k]}=(\zeta ,..,\zeta)\in \complex^{m-k}$ et
$\zeta=(x_{k+1}...x_{m})^{1/(m-k)}$.
\end{lem}
{\bf Preuve.} La d\émonstration se fait par r\écurrence, comme pour le lemme \ref{lem1}.
Supposons que pour $k+1\leq p < m$ et pour tout $(x_{1},..,x_{m})
\in\reel^{m}$ avec $x_{i}>0$ on ait
\begin{eqnarray}\label{eq22}
&&(\varphi -\tilde{\psi})([1,x_{1},..,x_{m}],[1,x_{1},..,x_{k}])\geq \nonumber \\
&& (\varphi -\tilde{\psi} )([1,x_{1},..,x_{k};(x_{k+1}...x_{p})^{\frac{1}{p-k}},..,
(x_{k+1}...x_{p})^{\frac{1}{p-k}},x_{p+1},..x_{m}],\nonumber\\
&&[1,x_{1},..,x_{k}]).
\end{eqnarray}
Cette propri\ét\é est claire pour $p=k+1$. Si 
l'in\égalit\é (\ref{eq22}) n'\était pas satisfaite au rang
$p+1$, il existerait alors un point
$(x_{1}^{0},..,x_{m}^{0})\in\reel^{m}$ avec $x_{i}^{0}>0$ pour
tout $i$, tel que
\begin{eqnarray}\label{eq23}
&&(\varphi -\tilde{\psi} )([1,x_{1}^{0},..,x_{m}^{0}],[1,x_{1}^{0},..,x_{k}^{0}])< \nonumber \\
&&(\varphi -\tilde{\psi})([1,x_{1}^{0},..,x_{k}^{0};(x_{k+1}^{0}...x_{p+1}^{0})^
{\frac{1}{p+1-k}},..,(x_{k+1}^{0}...x_{p+1}^{0})^{\frac{1}{p+1-k}},x_{p+2}^{0},..,x_{m}^{0}],\nonumber\\
&&[1,x_{1}^{0},..,x_{k}^{0}]).
\end{eqnarray}

En utilisant la continuit\é de $(\varphi -\tilde{\psi})$, on peut
supposer, quitte \`{a} en modifier l\ég\érement les
coordonn\ées, que le point $([1,x_{1}^{0},..,x_{m}^{0}],[1,x_{1}^{0},..,x_{k}^{0}])$ de
l'in\égalit\é (\ref{eq23}), v\érifie
$$(x_{1}^{0}..x_{k}^{0})^{1/(k+1)}\neq
(x_{k+1}^{0}..x_{m}^{0})^{1/(m-k)},$$ propri\ét\é dont on
aura besoin plus loin. En utilisant la $G$-invariance de
$\varphi$, on peut supposer que $x_{k+1}^{0}\leq ...\leq
x_{m}^{0}$. D'autre part, en tenant encore compte de la $G$
invariance de $\varphi$ et de l'hypoth\èse de r\écurrence
(\ref{eq22}) en les points
$$([1,x_{1}^{0},..,x_{k}^{0};x_{k+1}^{0},x_{k+2}^{0},..,x_{p}^{0},
x_{p+1}^{0},x_{p+2}^{0},..,x_{m}^{0}],[1,x_{1}^{0},..,x_{k}^{0}])$$ et
$$([1,x_{1}^{0},..,x_{k}^{0};x_{k+2}^{0},
x_{k+3}^{0},..,x_{p}^{0},x_{p+1}^{0},x_{k+1}^{0},x_{p+2}^{0},..,x_{m}^{0}],
[1,x_{1}^{0},..,x_{k}^{0}]),$$
on peut \écrire
\begin{eqnarray}\label{eq24}
&&(\varphi-\tilde{\psi} )([1,x_{1}^{0},..,x_{m}^{0}],[1,x_{1}^{0},..,x_{k}^{0}])\geq\nonumber\\
&& (\varphi-\tilde{\psi}
)([1,x_{1}^{0},..,x_{k}^{0};(x_{k+1}^{0}...x_{p}^{0})
^{\frac{1}{p-k}},..,(x_{k+1}^{0}...x_{p}^{0})^
{\frac{1}{p-k}},x_{p+1}^{0},x_{p+2}^{0},..,x_{m}^{0}],\nonumber\\
&&[1,x_{1}^{0},..,x_{k}^{0}])
\end{eqnarray}
et
\begin{eqnarray}\label{eq25}
&&(\varphi-\tilde{\psi})([1,x_{1}^{0},..,x_{k}^{0};x_{k+2}^{0},..,x_{p+1}^{0},x_{k+1}^{0},
x_{p+2}^{0},..,x_{m}^{0}],[1,x_{1}^{0},..,x_{k}^{0}])\geq \nonumber \\
&& (\varphi-\tilde{\psi})([1,x_{1}^{0},..,x_{k}^{0};(x_{k+2}^{0}..x_{p+1}^{0})
^{\frac{1}{p-k}},..,(x_{k+2}^{0}..x_{p+1}^{0})^{\frac{1}{p-k}},x_{k+1}^{0},
x_{p+2}^{0},..,x_{m}^{0}],\nonumber\\
&&[1,x_{1}^{0},..,x_{k}^{0}]).
\end{eqnarray}
Consid\érons maintenant la courbe $C$ d'\équation
$$t^{p-k}x=x_{k+1}^{0}...x_{p+1}^{0}$$ dans le plan r\éel
$\{[1,x_{1}^{0},..,x_{k}^{0},t,..,t,x,x_{p+2}^{0},..,x_{m}^{0}],[1,x_{1}^{0},..,x_{k}^{0}]\}$
param\étr\é par les variables $t$ et $x$. Les points
\begin{eqnarray}
\tilde{P}_{1}=([1,x_{1}^{0},..,x_{k}^{0};(x_{k+1}^{0}...x_{p}^{0})^
{\frac{1}{p-k}},..,(x_{k+1}^{0}...x_{p}^{0})^{\frac{1}{p-k}},
x_{p+1}^{0},x_{p+2}^{0},..,x_{m}^{0}],[1,x_{1}^{0},..,x_{k}^{0}])\nonumber
\end{eqnarray} 
et
\begin{eqnarray}
&&\tilde{P}_{2}=([1,x_{1}^{0},..,x_{k}^{0};(x_{k+2}^{0}...x_{p+1}^{0})
^{\frac{1}{p-k}},..,(x_{k+2}^{0}...x_{p+1}^{0})
^{\frac{1}{p-k}},x_{k+1}^{0},x_{p+2}^{0},..,x_{m}^{0}],\nonumber \\ &&[1,x_{1}^{0},..,x_{k}^{0}])\nonumber
\end{eqnarray}
appartiennent \`{a} la courbe $C$. Notons que les r\éels
$x_{i}^{0}$ pour $k+1\leq i\leq p+1$ ne sont pas tous \égaux, 
sinon (\ref{eq23}) deviendrait une \égalit\é.\\ Par suite,
sachant que l'on a choisi $x_{k+1}^{0}\leq ...\leq x_{p+1}^{0}$,
les points distincts $\tilde{P}_{1}$ et $\tilde{P}_{2}$ se trouvent strictement de
part et d'autre de la diagonale $t=x$ du plan pr\éc\édent.
\\Or la courbe $C$ coupe cette diagonale en le point

$$\tilde{P}_{3}=([1,x_{1}^{0},..,x_{k}^{0};
(x_{k+1}^{0}...x_{p+1}^{0})^{\frac{1}{p+1-k}},..,
(x_{k+1}^{0}...x_{p+1}^{0})^{\frac{1}{p+1-k}},x_{p+2}^{0},..,x_{m}^{0}],[1,x_{1}^{0},..,x_{k}^{0}])$$
qui intervient dans l'in\égalit\é (\ref{eq23}). D'autre part,
en utilisant les relations (\ref{eq23}), (\ref{eq24}) et (\ref{eq25})
on obtient :$$ (\varphi - \tilde{\psi})(\tilde{P}_{3})> (\varphi - \tilde{\psi}
)(\tilde{P}_{1})\mbox{ et }(\varphi - \tilde{\psi} )(\tilde{P}_{3})> (\varphi - \tilde{\psi}
)(\tilde{P}_{2}),$$ ce qui prouve que la fonction $(\varphi - \tilde{\psi} )$
admet un maximum local sur la courbe $C$. En cons\équence, la
restriction de la fonction $G$-invariante $(\varphi - \tilde{\psi} )$
\`{a} la courbe holomorphe (toujours not\'ee $C$) d'\équation
$\xi^{p-k}z=x_{k+1}^{0}...x_{p+1}^{0}$ du plan complexe
$\{([1,x_{1}^{0},..,x_{k}^{0};\xi,..,\xi,z,x_{p+2}^{0},..,x_{m}^{0}],
[1,x_{1}^{0},..,x_{k}^{0}])\}$
atteint un maximum local en un point $\tilde{P}=C(\zeta)$. Posons
$$C(\zeta)=([1,C^{1}(\zeta),..,C^{m}(\zeta)],[1,C^{1}(\zeta),..,C^{k}(\zeta)]),$$
$$\dot{C}^{\lambda}(\xi)=\frac{d C^{\lambda}}{d \xi}(\xi) \mbox{ et
} \dot{C}^{\overline{\mu}}(\xi)=\overline{\dot{C}^{\mu}(\xi)}.$$

Sachant que l'on a choisi le point $([1,x_{1}^{0},..,x_{m}^{0}],[1,x_{1}^{0},..,x_{k}^{0}])$ de
sorte que $$(x_{1}^{0}...x_{k}^{0})^{1/(k+1)}\neq
(x_{k+1}^{0}...x_{m}^{0})^{1/(m-k)},$$ l'\équation de la
courbe $C$ et les d\éfinitions de $\tilde{\psi}_{1}$ et $\tilde{\psi}_{2}$ montrent 
qu'en tout point de $C$ 
\begin{eqnarray} 
\tilde{\psi}_{1}([1,x_{1}^{0},..,x_{k}^{0};\xi,..,\xi,z,
x_{p+2}^{0},..,x_{m}^{0}],[1,x_{1}^{0},..,x_{k}^{0}])\neq\nonumber\\ 
\tilde{\psi}_{2}([1,x_{1}^{0},..,x_{k}^{0};\xi,..,\xi,z,x_{p+2}^{0},..,x_{m}^{0}],[1,x_{1}^{0},..,x_{k}^{0}]).
\end{eqnarray}
On peut alors supposer que $\tilde{\psi}=\tilde{\psi}_{1}$ dans un voisinage de $\tilde{P}$, 
la preuve \étant identique si l'on suppose $\tilde{\psi}=\tilde{\psi}_{2}$ dans
ce voisinage. On a donc :

$$\frac {\partial^{2}}{\partial \xi\partial \overline{\xi}}
\{(\varphi -\tilde{\psi}_{1})(C(\zeta))\} =\frac{\partial^{2}(\varphi
-\tilde{\psi}_{1})}{\partial z_{\lambda}
\partial \overline{z}_{\mu} }(C(\zeta))\dot{C}^{\lambda}(\zeta)
\dot{C}^{\overline{\mu}}(\zeta)$$ est n\égatif ou nul.
Comme $-\frac{\partial^{2}\tilde{\psi}_{1}}{\partial z_{\lambda}
\partial \overline{z}_{\mu}}= \tilde{g}_{\lambda\overline{\mu}},$
ceci exprime que la forme hermitienne de matrice:
$$(\tilde{g}_{\lambda\overline{\mu}}+\frac{\partial^{2}\varphi}{\partial
z_{\lambda}
\partial \overline{z}_{\mu}})_{\lambda ,\mu}=(\frac{\partial^{2}(\varphi -\tilde{\psi}_{1} )}
{\partial z_{\lambda}
\partial \overline{z}_{\mu}})_{\lambda,\mu}$$
est n\égative en $\tilde{P} = C(\zeta)$. On en
d\éduit une contradiction avec la $\tilde{g}$-admissibilit\é de
$\varphi$ en $\tilde{P}$. D'o\`u l'in\'egalit\'e (\ref{eq22}) au rang $p+1$ et, par
cons\'equent, le lemme \ref{lem21}. 


\begin{lem}\label{lem22}
Etant donn\ée une fonction $\varphi\in C^{\infty}(X)$,
$\tilde{g}$-admissible, G-invariante. Si $x_{i}=\mid z_{i}\mid >0$ pour
tout $i$, on a:
\begin{eqnarray}\label{eq26}
(\varphi -\tilde{\psi}
)([x_{0},x_{1},..,x_{k};1,x_{k+2},..,x_{m}],
[x_{0},x_{1},..,x_{k}])\geq\nonumber\\ (\varphi -\tilde{\psi}
)([\eta,\eta,..,\eta;1,x_{k+2},..,x_{m}],[1^{[k+1]}]),
\end{eqnarray}
o\`u $\eta=(x_{0}x_{1}...x_{k})^{1/(k+1)}$.
\end{lem}
{\bf Preuve.} Comme dans le lemme \ref{lem21}, la preuve s'effectue
par r\écurrence. Supposons que pour $0\leq p < k$ et
pour tout $(x_{0},..,x_{k};x_{k+2},..,x_{m}) \in\reel^{m}$ avec
$x_{i}>0$, on ait
\begin{eqnarray}\label{eq27}
&&(\varphi -\tilde{\psi} )([x_{0},..,x_{k},1,x_{k+2},..,x_{m}],[x_{0},..,x_{k}])\geq
\nonumber\\
&&(\varphi -\tilde{\psi}
)([(x_{0}...x_{p})^{\frac{1}{p+1}},..,
(x_{0}...x_{p})^{\frac{1}{p+1}},x_{p+1},..,x_{k};1,x_{k+2},..,x_{m}],\nonumber\\
&&[(x_{0}...x_{p})^{\frac{1}{p+1}},..,
(x_{0}...x_{p})^{\frac{1}{p+1}},x_{p+1},..,x_{k}]).
\end{eqnarray}
Cette hypoth\èse est v\érifi\ée pour $p=0$. Si
l'in\égalit\é (\ref{eq27}) n'\était pas satisfaite au rang
$p+1$, il existerait un point
$(x_{0}^{0},..,x_{k};x_{k+2},..,x_{m}^{0})\in\reel^{m}$ avec
$x_{i}^{0}>0$ pour tout $i$, tel que :
\begin{eqnarray}\label{eq28}
&&(\varphi -\tilde{\psi}
)([x_{0}^{0},..,x_{k}^{0};1,x_{k+2}^{0},..,x_{m}^{0}],[x_{0}^{0},..,x_{k}^{0}])< \nonumber\\
&& (\varphi
-\tilde{\psi})([(x_{0}^{0}..x_{p+1}^{0})^{\frac{1}{p+2}},..,(x_{0}^{0}..x_{p+1}^{0})^{\frac{1}{p+2}},
x_{p+2}^{0},...,x_{k}^{0};1,x_{k+2}^{0},..,x_{m}^{0}],\nonumber\\
&&[(x_{0}^{0}..x_{p+1}^{0})^{\frac{1}{p+2}},..,(x_{0}^{0}..x_{p+1}^{0})^{\frac{1}{p+2}},
x_{p+2}^{0},...,x_{k}^{0}]).
\end{eqnarray}
Comme au lemme \ref{lem21}, on peut supposer que le point
$$([x_{0}^{0},..,x_{k}^{0};1,x_{k+2}^{0},..,x_{m}^{0}],[x_{0}^{0},..,x_{k}^{0}])\in X $$ v\érifie
$$(x_{0}^{0}..x_{k}^{0})^{1/(k+1)}\neq
(x_{k+2}^{0}..x_{m}^{0})^{1/(m-k)},$$ et que $x_{0}^{0}\leq
...\leq x_{p+1}^{0}$. D'autre part, en tenant compte de la $G$-invariance 
de $\varphi$ et de l'hypoth\èse de r\écurrence
(\ref{eq27}) en les points 
$$([x_{0}^{0},x_{1}^{0},..,x_{p}^{0},x_{p+1}^{0},..,x_{k}^{0};1,x_{k+2}^{0},..,x_{m}^{0}],
[x_{0}^{0},x_{1}^{0},..,x_{p}^{0},x_{p+1}^{0},..,x_{k}^{0}])$$
et
$$([x_{1}^{0},..,x_{p+1}^{0},x_{0}^{0},..,x_{k}^{0};1,x_{k+2}^{0},..,x_{m}^{0}],
[x_{1}^{0},..,x_{p+1}^{0},x_{0}^{0},..,x_{k}^{0}]),$$
on a 
\begin{eqnarray}\label{eq29}
&&(\varphi-\tilde{\psi} )([x_{0}^{0},..,x_{k}^{0};1,x_{k+2}^{0},
..,x_{m}^{0}],[x_{0}^{0},..,x_{k}^{0}])\geq\nonumber\\
&& (\varphi-\tilde{\psi}
)([(x_{0}^{0}..x_{p}^{0})^{\frac{1}{p+1}},..,(x_{0}^{0}..x_{p}^{0})^{\frac{1}{p+1}},
x_{p+1}^{0},..,x_{k}^{0};1,x_{k+2}^{0},..,x_{m}^{0}],\nonumber\\
&&[(x_{0}^{0}..x_{p}^{0})^{\frac{1}{p+1}},..,(x_{0}^{0}..x_{p}^{0})^{\frac{1}{p+1}},
x_{p+1}^{0},..,x_{k}^{0}]) 
\end{eqnarray}
et
\begin{eqnarray}\label{eq210}
&&(\varphi-\tilde{\psi})([x_{1}^{0},..,x_{p+1}^{0},x_{0}^{0},..,x_{k}^{0};1,x_{k+2}^{0},..,x_{m}^{0}],
[x_{1}^{0},..,x_{p+1}^{0},x_{0}^{0},..,x_{k}^{0}])\geq
\nonumber \\
&&(\varphi-\tilde{\psi})([(x_{1}^{0}..x_{p+1}^{0})^{\frac{1}{p+1}},..,(x_{1}^{0}..x_{p+1}^{0})^
{\frac{1}{p+1}},
x_{0}^{0},x_{p+2}^{0},..,x_{k}^{0};1,x_{k+2}^{0},..,x_{m}^{0}],\nonumber\\
&&[(x_{1}^{0}..x_{p+1}^{0})^{\frac{1}{p+1}},..,(x_{1}^{0}..x_{p+1}^{0})^
{\frac{1}{p+1}},
x_{0}^{0},x_{p+2}^{0},..,x_{k}^{0}]).
\end{eqnarray}
Consid\érons maintenant la courbe $C$ d'\équation
$$t^{p+1}x=x_{0}^{0}...x_{p+1}^{0}$$ du plan r\éel
$\{([t,..,t,x,x_{p+2}^{0},..,x_{k}^{0};1,x_{k+2}^{0},..,x_{m}^{0}],[t,..,t,x,x_{p+2}^{0},..,x_{k}^{0}])\}$
param\étr\é par les variables $t$ et $x$. Les points 
\begin{eqnarray}
&&\tilde{Q}_{1}=([(x_{0}^{0}..x_{p}^{0})^{\frac{1}{p+1}},..,
(x_{0}^{0}..x_{p}^{0})^{\frac{1}{p+1}},x_{p+1}^{0},x_{p+2}^{0},..,x_{k}^{0};1,
x_{k+2}^{0},..,x_{m}^{0}],\nonumber\\
&&[(x_{0}^{0}..x_{p}^{0})^{\frac{1}{p+1}},..,
(x_{0}^{0}..x_{p}^{0})^{\frac{1}{p+1}},x_{p+1}^{0},x_{p+2}^{0},..,x_{k}^{0}])\nonumber
\end{eqnarray} 
et
\begin{eqnarray}
&&\tilde{Q}_{2}=([(x_{1}^{0}..x_{p+1}^{0})^{\frac{1}{p+1}},..,(x_{1}^{0}..x_{p+1}^{0})^{\frac{1}{p+1}},
x_{0}^{0},x_{p+2}^{0},..,x_{k}^{0};1,x_{k+2}^{0},..,x_{m}^{0}],\nonumber\\
&&[(x_{1}^{0}..x_{p+1}^{0})^{\frac{1}{p+1}},
..,(x_{1}^{0}..x_{p+1}^{0})^{\frac{1}{p+1}},
x_{0}^{0},x_{p+2}^{0},..,x_{k}^{0}])\nonumber
\end{eqnarray} 
appartiennent \`{a} la courbe $C$.\\D'autre part 
 les r\éels $x_{i}^{0}$ pour $0\leq i\leq p+1$ ne sont pas tous
\égaux, sinon (\ref{eq28}) serait une \égalit\é.\\
Par suite, sachant que l'on a choisi $x_{0}^{0}\leq ...\leq
x_{p+1}^{0}$, les points distincts $\tilde{Q}_{1}$ et $\tilde{Q}_{2}$ se trouvent strictement 
de part et
d'autre de la diagonale $t=x$ du plan pr\éc\édent. \\Or la
courbe $C$ coupe cette diagonale en le point 

\begin{eqnarray}
&&\tilde{Q}_{3}=([(x_{0}^{0}...x_{p+1}^{0})^{\frac{1}{p+2}},..,
(x_{0}^{0}...x_{p+1}^{0})^{\frac{1}{p+2}},x_{p+2}^{0},..,x_{k}^{0};
1,x_{k+2}^{0},..,x_{m}^{0}],\nonumber\\
&&[(x_{0}^{0}...x_{p+1}^{0})^{\frac{1}{p+2}},..,
(x_{0}^{0}...x_{p+1}^{0})^{\frac{1}{p+2}},x_{p+2}^{0},..,x_{k}^{0}])\nonumber
\end{eqnarray} qui intervient dans
l'in\égalit\é (\ref{eq28}). D'autre part, les relation
(\ref{eq28}), (\ref{eq29}) et (\ref{eq210}) donnent 
$$ (\varphi -
\tilde{\psi})(\tilde{Q}_{3})> (\varphi - \tilde{\psi} )(\tilde{Q}_{1})\mbox{ et }(\varphi - \tilde{\psi}
)(\tilde{Q}_{3})> (\varphi - \tilde{\psi} )(\tilde{Q}_{2}),$$ ce qui prouve que la
fonction $(\varphi - \tilde{\psi} )$ admet un maximum local sur la courbe
$C$.\\
Sachant que l'on a choisi le point
$$([x_{0}^{0},..,x_{k}^{0};1,x_{k+2}^{0},..,x_{m}^{0}],[x_{0}^{0},..,x_{k}^{0}])$$ de sorte que
$$(x_{0}^{0}..x_{k}^{0})^{1/(k+1)}\neq
(x_{k+2}^{0}..x_{m}^{0})^{1/(m-k)},$$
on conclut de la
m\^{e}me mani\ère qu'au lemme pr\éc\édent en consid\érant la restriction de 
$(\varphi -\tilde{\psi})$ \`{a} une courbe holomorphe convenable.


A l'instar du lemme \ref{coro4}, les lemmes \ref{lem21} et \ref{lem22} permettent d'\établir le
\begin{lem}\label{coro24}
Etant donn\ée une fonction $\varphi\in C^{\infty}(X)$,
$\tilde{g}$-admissible, G-invariante, avec $x_{i}=\mid z_{i}\mid >0$ pour
tout $i$, on a :
\begin{eqnarray}\label{eq211}
&&(\varphi -\tilde{\psi}
)([1,x_{1},..,x_{k};x_{k+1},x_{k+2},..,x_{m}],[1,x_{1},..,x_{k}])\nonumber\\
&&\geq(\varphi -\tilde{\psi} )([1^{[k+1]};\nu^{[m-k]}],[1^{[k+1]}]),
\end{eqnarray}
o\`u
$\nu=\frac{(x_{k+1}...x_{m})^{1/(m-k)}}{(x_{1}...x_{k})^{1/(k+1)}}$.
\end{lem}
Montrons maintenant le 

\begin{lem}\label{lem24}
Etant donn\ée une fonction $\varphi\in C^{\infty}(X)$,
$\tilde{g}$-admissible, G-invariante, $\forall \zeta >0$, on a :
\begin{eqnarray}\label{eq212}
(\varphi -\tilde{\psi} )([1^{[k+1]};\zeta^{[m-k]}],[1^{[k+1]}])\geq 0,
\end{eqnarray}
\end{lem}


{\bf Preuve.} On raisonne sur la position du point $\tilde{R}_{0}\in\proj_{m}\complex$
o\`{u} $\varphi$ atteint son maximum. En vertu de la $G$-invariance
de $\varphi$, on peut supposer qu'il s'\écrit sous la forme 
$$\tilde{R}_{0}=([y^{0}_{0},..,y^{0}_{k};y^{0}_{k+1},..,y^{0}_{m}],[\rho ^{0}_{0},..,\rho ^{0}_{k}]),$$
o\`{u} les $y^{0}_{i}$ et $\rho ^{0}_{i}$ sont des r\éels positifs v\érifiant 
$y^{0}_{0}\geq y^{0}_{1}\geq..\geq y^{0}_{k}$, $y^{0}_{k+1}\geq
y^{0}_{k+2}\geq..\geq y^{0}_{m}$, $\rho ^{0}_{0}\geq..\geq\rho ^{0}_{k}$ et o\ù $(\rho ^{0}_{0},..,\rho ^{0}_{k})$ et 
$(y^{0}_{0},..,y^{0}_{k})$ sont 
parall\èles. Deux cas se pr\ésentent : ou
bien l'un des $y^{0}_{0},..,y^{0}_{k}$ est non nul, ou bien tous
les $y^{0}_{0},..,y^{0}_{k}$ sont nuls.
\\{\bf{Cas A}} :
l'un des $y^{0}_{0},..,y^{0}_{k}$ est non nul. On peut alors se
placer dans la carte de $X$ d\écrite par les points 
$$([z_{0},..,z_{k},z_{k+1},..,z_{m}],[\xi_{0},..,\xi_{k}])\in
\proj_{m}\complex \times\proj_{k}\complex$$ 
tels que $z_{0}\neq 0$.
Ceci permet d'\écrire le point
$\tilde{R}_{0}$ sous la forme 
$$\tilde{R}_{0}=([1,x^{0}_{1},..,x^{0}_{k};x^{0}_{k+1},..,x^{0}_{m}],[1,x^{0}_{1},..,x^{0}_{k}]),$$
o\`{u} les r\éels positifs $x^{0}_{i}$ v\érifient :
$1\geq x^{0}_{1}\geq..\geq x^{0}_{k}$ et $x^{0}_{k+1}\geq ...\geq
x^{0}_{m}$.
Raisonnons par l'absurde, et supposons qu'il existe un point
$$\tilde{R}_{1}=([1^{[k+1]};\zeta_{0}^{[m-k]}],[1^{[k+1]}])$$ tel que
l'on ait $\zeta_{0}>0$ et 
\begin{equation}\label{p21}
(\varphi -\tilde{\psi} )(\tilde{R}_{1})<0.
\end{equation}
On envisage alors les deux sous-cas suivants : $x^{0}_{k+1}<\zeta_{0}$ 
puis $x^{0}_{k+1}\geq \zeta_{0}$.
\begin{itemize}
    \item \underline{$x^{0}_{k+1}<\zeta_{0}$}.
\end{itemize}
On introduit alors la fonction auxiliaire
$$\tilde{\psi}_{0}=\log\frac{\mid z_{0}\mid^{2(m+1-k)}\mid \xi_{0}\mid^{2k}}{
    (\mid z_{0}\mid^{2}+...+\mid z_{m}\mid^{2})^{m+1-k}(\mid \xi_{0}\mid^{2}+...+\mid \xi_{k}\mid^{2})^{k}}.$$
D'une part, puisque $\varphi\leq 0$,
\begin{equation}\label{p20}
(\varphi -\tilde{\psi}_{0})([1,0^{[m]}],[1,0^{[k]}])=\varphi 
([1,0^{[m]}],[1,0^{[k]}])\leq 0.
\end{equation}
De plus, sachant que $\varphi (\tilde{R}_{0})=0$ et $\tilde{\psi}_{0}\leq 0$, on a
\begin{equation}\label{p23}
(\varphi -\tilde{\psi}_{0})(\tilde{R}_{0})\geq 0.
\end{equation}
Si $\tilde{R}_{0}\neq ([1,0^{[m]}],[1,0^{[k]}])$, $\tilde{\psi}_{0}(\tilde{R}_{0})< 0$ et l'in\égalit\é (\ref{p23}) est
alors stricte. Si $\tilde{R}_{0} = ([1,0^{[m]}],[1,0^{[k]}])$, 
quitte \`{a} se placer en un point $\tilde{R}$ arbitrairement voisin de
$\tilde{R}_{0}$, on peut supposer $(\varphi -\tilde{\psi}_{0})(\tilde{R})> 0$. En effet, si dans
un voisinage de $\tilde{R}_{0}$ on avait $(\varphi -\tilde{\psi}_{0})\leq 0$, comme 
$(\varphi -\tilde{\psi}_{0})(\tilde{R}_{0})= 0$, $(\varphi -\tilde{\psi}_{0})$ admettrait
alors un maximun local en $\tilde{R}_{0}$, ce qui mettrait en
d\éfaut l'admissibilit\é de $\varphi$ en ce point, sachant que 
$$\partial_{\lambda\overline{\mu}}(\varphi -\tilde{\psi}_{0})(\tilde{R}_{0})=
(\tilde{g}_{\lambda\overline{\mu}}+\partial_{\lambda\overline{\mu}}\varphi )(\tilde{R}_{0}).$$ Dans tous les cas, 
on peut
donc affirmer qu'il existe un point $$\tilde{R}'_{0}=([1,{a}_{1},..,{a}_{m}],[1,{a}_{1},..,{a}_{k}])$$
v\érifiant
\begin{equation}\label{p23'}
(\varphi -\tilde{\psi}_{0})(\tilde{R}'_{0})> 0.
\end{equation}
Par continuit\é et $G$-invariance de $\varphi$, on peut supposer $1>a_{1}\geq...\geq a_{k}>0$ et 
$\zeta_{0}> a_{k+1}\geq...\geq a_{m}>0$.
D'autre part, l'in\égalit\é (\ref{p21}) jointe aux d\'efinitions de $\tilde{R}_1$, $\tilde{\psi}_0$, $\tilde{\psi}_1$ 
et $\tilde{\psi}=\inf (\tilde{\psi}_1 , \tilde{\psi}_2)$ implique
\begin{equation}\label{p21'}
(\varphi -\tilde{\psi}_{0})(\tilde{R}_{1})=(\varphi -\tilde{\psi}_{1})(\tilde{R}_{1})\leq(\varphi
-\tilde{\psi})(\tilde{R}_{1})< 0.
\end{equation}

La courbe :
\begin{eqnarray}
&[0,1]\ni t\rightarrow ([1,t,t^{(\ln a_{2})/(\ln a_{1})},..,t^{(\ln a_{k})/(\ln a_{1})};
\zeta_{0}t^{\frac{\ln(a_{k+1}/\zeta_{0})}{\ln
a_{1}}},..,\zeta_{0}t^{\frac{\ln(a_{m}/\zeta_{0})}{\ln
a_{1}}}],&\nonumber\\
&[1,t,t^{(\ln a_{2})/(\ln a_{1})},..,t^{(\ln a_{k})/(\ln a_{1})}])&
\nonumber
\end{eqnarray} 
passe par $([1,0^{[m]}],[1,0^{[k]}])$ en $t=0$ puis par ${\tilde{R}'}_{0}$
en $t=a_{1}$ et enfin par le point $\tilde{R}_{1}$ en $t=1$, valeurs
en lesquelles, d'apr\ès (\ref{p20}), (\ref{p23'}) et (\ref{p21'}),
$(\varphi-\tilde{\psi}_{0})$ est respectivement n\égative, positive
puis \`{a} nouveau n\égative. L'invariance de cette fonction
par l'action des $\exp(i\theta)$, permet donc de d\éduire que
$(\varphi-\tilde{\psi}_{0})$ atteint un maximum sur la courbe holomorphe,
complexifi\ée de la courbe d\écrite plus haut, ce qui
contredit encore une fois l'admissibilit\é de $\varphi$.
\begin{itemize}
\item \underline{$x^{0}_{k+1}\geq\zeta_{0}$}.
\end{itemize}
D\ésignons dans ce cas par $p\in \{1,..,m-k\}$ l'entier pour
lequel on a $$x^{0}_{k+1}\geq ...\geq x^{0}_{k+p}>\zeta_{0} \mbox{
et } \zeta_{0}\geq x^{0}_{k+p+1}\geq ...\geq x^{0}_{m},$$ et
consid\érons la fonction auxiliaire

$$\tilde{\psi}_{k+1}=\log\frac{\mid z_{k+1}\mid^{2(m+1-k)}\mid\xi_{0}\mid^{2k}}{
(\mid z_{0}\mid^{2}+...+\mid z_{m}\mid^{2})^{m+1-k}(\mid \xi_{0}\mid^{2}+...+\mid \xi_{k}\mid^{2})^{k}}.$$
On a
\begin{equation}\label{p23''}
(\varphi -\tilde{\psi}_{k+1})(\tilde{R}_{0})> 0.
\end{equation}
La fonction $(\varphi -\tilde{\psi}_{k+1})$ \étant continue, quitte
\`{a} se placer en un point voisin de $\tilde{R}_{0}$, on peut supposer
tous les $x_{i}^{0}$ non nuls. Posons alors:
$$\alpha_{2}=\frac{\ln x_{2}^{0}}{\ln x_{1}^{0}},..,\alpha_{k}=\frac{\ln x_{k}^{0}}{\ln x_{1}^{0}};
\alpha_{k+1}=\frac{\ln(x_{k+1}^{0}/\zeta_{0})}{\ln x_{1}^{0}},..,
\alpha_{m}=\frac{\ln(x_{m}^{0}/\zeta_{0})}{\ln x_{1}^{0}}.$$
Sachant que l'on a $1\geq x_{1}^{0}\geq..\geq x_{k}^{0}$;
$x_{k+1}^{0}\geq..\geq x_{k+p}^{0}\geq\zeta_{0}$ et $\zeta_{0}\geq
x_{k+p+1}^{0}\geq..\geq x_{m}^{0}$, on en d\éduit que 
$\alpha_{2},..,\alpha_{k}\geq 0$, $\alpha_{k+1}\leq ...\leq
\alpha_{p+k}\leq 0$ et $\alpha_{p+k+1},..,\alpha_{m}\geq 0$, d'o\`{u}, en notant 
$$\tilde{R}_{\varepsilon}=([1,\varepsilon,\varepsilon^{\alpha_{2}},..,\varepsilon^{\alpha_{k}},
\zeta_{0}\varepsilon^{\alpha_{k+1}},..,\zeta_{0}\varepsilon^{\alpha_{m}}],
[1,\varepsilon,\varepsilon^{\alpha_{2}},..,\varepsilon^{\alpha_{k}}]),$$ on a 
\begin{eqnarray*}
\lim_{\varepsilon\rightarrow
0}\tilde{\psi}_{k+1}(\tilde{R}_{\varepsilon }) 
&=& \lim_{\varepsilon\rightarrow
0}\ln\frac{\zeta_{0}^{2(m+1-k)}\varepsilon^{2(m+1-k)\alpha_{k}}}
{[1+\varepsilon^{2}+\varepsilon^{2\alpha_{2}}+..+\varepsilon^{2\alpha_{k}}+
\zeta_{0}^{2}(\varepsilon^{2\alpha_{k+1}}+..+\varepsilon^{2\alpha_{m}})]^{m+1-k}}  \\
&=&\ln\lim_{t\rightarrow
\infty}\frac{\zeta_{0}^{2(m+1-k)}t^{-2\alpha_{k+1}(m+1-k)}}
{[\zeta_{0}^{2}(t^{-2\alpha_{k+1}}+t^{-2\alpha_{k+2}}+..+t^{-2\alpha_{p}})]^{(m+1-k)}}
= \ln 1=0,
\end{eqnarray*} 
$-\alpha_{k+1}$ \étant la plus grande puissance intervenant au
d\énominateur. Sachant que $\varphi ([\tilde{R}_{\varepsilon}])\leq 0$ et compte tenu de (\ref{p23''}) il
existe $\varepsilon_{0}$ tel que l'on ait
\begin{equation}\label{p24}
(\varphi
-\tilde{\psi}_{k+1})(\tilde{R}_{\varepsilon_{0} })\leq -\tilde{\psi}_{k+1}(\tilde{R}_{\varepsilon_{0} })<
(\varphi -\tilde{\psi}_{k+1})(\tilde{R}_{0}).
\end{equation}
D'autre part, l'in\'egalit\'e (\ref{p21}), jointe aux d\'efinitions de $\tilde{R}_1$, $\tilde{\psi}_{k+1}$, $\tilde{\psi}_{2}$ 
et $\tilde{\psi} = \inf (\tilde{\psi}_{1},\tilde{\psi}_{2})$ donne~: 
\begin{equation}\label{p25}
(\varphi
-\tilde{\psi}_{k+1})(\tilde{R}_{1})= (\varphi-\tilde{\psi}_{2})(\tilde{R}_{1})\leq
(\varphi -\tilde{\psi})(\tilde{R}_{1})<0.
\end{equation}
La courbe
$$[\varepsilon_{0},1]\ni t\rightarrow ([1,t,t^{\alpha_{2}},..,t^{\alpha_{k}},
\zeta_{0}t^{\alpha_{k+1}},..,\zeta_{0}t^{\alpha_{m}}],
[1,t,t^{\alpha_{2}},..,t^{\alpha_{k}}]),$$ passe par
$\tilde{R}_{\varepsilon_{0}}$
en $t=\varepsilon_{0}$ puis par $\tilde{R}_{0}$ en $t=x_{1}^{0}$ et enfin
par $\tilde{R}_{1}$ en $t=1$, ce qui, en vertu de 
(\ref{p24}), (\ref{p23''}) et (\ref{p25}) prouve l'existence d'un maximum local pour la 
fonction $(\varphi-\tilde{\psi}_{k+1})$ sur la
courbe pr\écit\ée. Ceci contredit, \`{a} l'instar du cas
pr\éc\édent, l'hypoth\èse d'admissibilit\é de la
fonction $\varphi$.

{\bf{Cas B :}}  $y^{0}_{0}=...=y^{0}_{k}=0$. On peut alors se
placer dans la carte de $X$ d\écrite par les points 
$$([z_{0},..,z_{k};z_{k+1},..,z_{m}],[\xi_{0},..,\xi_{k}])\in
\proj_{m}\complex \times\proj_{k}\complex$$ 
tels que $z_{k+1}\neq 0$ et $\xi_{0}\neq 0$, de sorte que
le point $\tilde{R}_{0}$ o\`u $\varphi$ atteint son maximum \'egal \`a z\'ero puisse 
s'\'ecrire sous la forme 
$$\tilde{R}_{0}=([0,0,..,0;1,x_{k+2}^{0},..,x_{m}^{0}],[1,u_{1}^{0},..,u_{k}^{0}]).$$
On peut aussi supposer, en utilisant la $G$-invariance de $\varphi$, que 
$1\geq x_{k+2}^{0}\geq..\geq x_{m}^{0}$ et $1\geq u_{1}^{0}\geq..\geq u_{k}^{0}$. 
On montrera une version \'equivalente du lemme \ref{lem24}, \`a savoir que
\begin{equation}\label{equ20}
(\varphi -\tilde{\psi} )([\zeta^{[k+1]};1^{[m-k]}],[1^{[k+1]}])\geq 0
\end{equation}
pour tout $\zeta >0$. 

Raisonnons par l'absurde, et supposons qu'il existe un point
$$\tilde{R}_{k+1}=([\zeta_{0}^{[k+1]};1^{[m-k]}],[1^{[k+1]}])$$ de $X$ tel que
l'on ait $\zeta_{0}>0$ et 
\begin{equation}\label{equ21}
(\varphi -\tilde{\psi} )(\tilde{R}_{k+1})<0.
\end{equation}
On consid\`ere alors la fonction auxiliaire
$\tilde{\psi}_{k+1}$ introduite plus haut.
Sachant que $\varphi\leq 0$, on a 
\begin{equation}\label{equ22}
(\varphi -\tilde{\psi}_{k+1})([0^{[k+1]};1,0,..,0],[1,0^{[k]}])=\varphi ([0^{[k+1]};1,0,..,0],[1,0^{[k]}])\leq 0.
\end{equation}
D'autre part, sachant que $\varphi (\tilde{R}_{0})=0$ et $\tilde{\psi}_{k+1}\leq 0$, on a

\begin{equation}
(\varphi -\tilde{\psi}_{k+1})(\tilde{R}_{0})=-\tilde{\psi}_{k+1}(\tilde{R}_{0})\geq 0.
\end{equation}
Cette in\'egalit\'e est stricte d\`es que $$\tilde{R}_{0}\neq ([0^{[k+1]};1,0,..,0],[1,0^{[k]}]).$$ Si 
$\tilde{R}_{0} = ([0^{[k+1]};1,0,..,0],[1,0^{[k]}])$, quitte \`a se placer en un point arbitrairement voisin de 
$\tilde{R}_{0}$, on peut supposer la derni\`ere in\'egalit\'e stricte. En effet, si dans un voisinage de $\tilde{R}_{0}$, on avait 
$\varphi -\tilde{\psi}_{k+1}\leq 0$, alors $\varphi -\tilde{\psi}_{k+1}$ admettrait un maximum local en  
$\tilde{R}_{0}$, ce qui contredirait l'admissibilit\'e de $\varphi$ en $\tilde{R}_{0}$. A l'instar du cas A, 
il existe donc un point  $$\tilde{R}'_{0}=([{c}_{0},..,{c}_{k};1,{c}_{k+2},..,{c}_{m}],[{c}_{0},..,{c}_{k}])$$
v\érifiant
\begin{equation}\label{equ23}
(\varphi -\tilde{\psi}_{k+1})(\tilde{R}'_{0})> 0.
\end{equation}
Par continuit\é et $G$-invariance de $\varphi$, on peut supposer 
$\zeta_{0}> c_{0}>...> c_{k}>0$ et $1>c_{k+2}>...> c_{m}>0$. 
D'autre part, l'in\égalit\é (\ref{equ21}) jointe aux d\'efinitions de $\tilde{R}_{k+1}$, $\tilde{\psi}_{k+1}$, $\tilde{\psi}_2$ 
et $\tilde{\psi}=\inf (\tilde{\psi}_1 , \tilde{\psi}_2)$ implique
\begin{equation}\label{equ24}
(\varphi -\tilde{\psi}_{k+1})(\tilde{R}_{k+1})=(\varphi -\tilde{\psi}_{2})(\tilde{R}_{k+1})\leq(\varphi
-\tilde{\psi})(\tilde{R}_{k+1})< 0.
\end{equation}

La courbe de $X$~:
\begin{eqnarray}
&[0,1]\ni t\rightarrow ([\zeta_{0}t^{\frac{\ln(c_{0}/\zeta_{0})}{\ln
c_{k+2}}},..,\zeta_{0}t^{\frac{\ln(c_{k}/\zeta_{0})}{\ln
c_{k+2}}};1,t,t^{(\ln c_{k+3})/(\ln c_{k+2})},..,t^{(\ln c_{m})/(\ln c_{k+2})}],&\nonumber\\
&[1,t^{\frac{\ln(c_{1}/c_{0})}{\ln
c_{k+2}}},..,t^{\frac{\ln(c_{k}/c_{0})}{\ln
c_{k+2}}}])&\nonumber
\end{eqnarray} 
passe par $([0^{[k+1]},1,0,..,0],[1,0^{[k]}])$ en $t=0$ puis par $\tilde{R}_{0}$
en $t=c_{k+2}$ et enfin par le point $\tilde{R}_{k+1}$ en $t=1$, valeurs
en lesquelles, d'apr\ès (\ref{equ22}), (\ref{equ23}) et (\ref{equ24}),
$(\varphi-\tilde{\psi}_{k+1})$ est respectivement n\égative, positive
puis n\égative. L'invariance de cette fonction
par l'action des $\exp(i\theta)$, permet donc de d\éduire que
$(\varphi-\tilde{\psi}_{0})$ atteint un maximum sur la courbe holomorphe,
d\'eduite de la courbe d\écrite plus haut, ce qui
contredit l'admissibilit\é de $\varphi$, d'o\`u (\ref{equ20}) 
et le lemme \ref{lem24}.


{\subsection{Preuve du corollaire \ref{coro2}.}}

Soit $\varphi\in C^{\infty}(X)$ une fonction
$\tilde{g}$-admissible et $G$-invariante, dont le $\sup$ sur
$X$ est nul. D'apr\`es le th\'eor\`eme \ref{th2},
on a $\varphi \geq \tilde{\psi}$ et par suite, pour tout $\alpha \geq 0$,
$$\int_{X} \exp(-\alpha\varphi ) dv \leq
\int_{X}\exp(-\alpha\tilde{\psi} )dv .$$ Cherchons les
valeurs de $\alpha$ pour lesquelles cette derni\`ere int\'egrale
converge. Pour ce faire, on estimera
$\int_{X}\exp(-\alpha\tilde{\psi}_{1} )dv$ et $\int_{X}\exp(-\alpha\tilde{\psi}_{2} )dv$ 
dans la carte
dense correspondant \à la param\étrisation 
$$([1,z_{1},..,z_{m}],[1,z_{1},..,z_{k}]).$$ Dans cette carte,
l'\'el\'ement de volume est donn\'e par (c.f. \cite{BC})~: 
$$dv=\det((\tilde{g}_{\lambda \overline{\mu}}))dz_{1}\wedge d\overline{z}_{1}\wedge ...\wedge
dz_{m}\wedge d\overline{z}_{m},$$ o\ù
\begin{eqnarray}
&&\det((\tilde{g}_{\lambda \overline{\mu}}))=(-1)^{m}\frac{(m+1-k)^{m-k}} {(1+\mid z_{1}\mid^{2}+...+\mid
z_{k}\mid^{2})^{k}(1+\mid z_{1}\mid^{2}+...+\mid
z_{m}\mid^{2})^{m+1}}\nonumber\\
&&\times [k(1+\mid z_{1}\mid^{2}+...+\mid
z_{m}\mid^{2})+(m-k+1)(1+\mid z_{1}\mid^{2}+...+\mid
z_{k}\mid^{2})]^{k}.\nonumber
\end{eqnarray} 
En utilisant le
fait que $\tilde{\psi}_{1}$ et $\tilde{\psi}_{2}$ ne d\'ependent que des $\mid
z_{p}\mid$, le changement de variables $u_{p}=\mid z_{p}\mid^{2}$
donne~:
\begin{eqnarray}
&&\int_{X}\exp(-\alpha\tilde{\psi}_{1} )dv= Cst\sum_{i=0}^{k}C^{i}_{k}\times\nonumber\\
&&\int_{0}^{+\infty}..\int_{0}^{+\infty}\frac{
(u_{1}...u_{k})^{-\alpha(m+1)/(k+1)}du_{1}...du_{m}} {(1+
u_{1}+...+u_{m})^{(1-\alpha)(m+1-k)+i}
(1+u_{1}+..+u_{k})^{(1-\alpha)k-i}}\nonumber\\
&&= \int_{0}^{+\infty}...
\int_{0}^{+\infty}\frac{
(u_{1}...u_{k})^{-\alpha(m+1)/(k+1)}du_{1}...du_{k}} {(1+
u_{1}+...+u_{m})^{(1-\alpha)(m+1)-(m-k)}}.\nonumber
\end{eqnarray}
qui converge, ind\épendemment de $i$ pour $\alpha < (k+1)/(m+1)$. D'autre part, 
\begin{eqnarray}
&&\int_{X}\exp(-\alpha\tilde{\psi}_{2} )dv= Cst\sum_{i=0}^{k}C^{i}_{k}\times\nonumber\\
&&\int_{0}^{+\infty}..\int_{0}^{+\infty}\frac{
(u_{1}...u_{k})^{-\alpha k/(k+1)}(u_{k+1}...u_{m})^{-\alpha(m+1-k)/(m-k)}du_{1}...du_{m}} {(1+
u_{1}+...+u_{m})^{(1-\alpha)(m+1-k)+i}
(1+u_{1}..+u_{k})^{(1-\alpha)k-i}}.\nonumber
\end{eqnarray}
Pour cette derni\ère int\égrale, la convergence en z\éro exige la condition 
$$\alpha <\inf\{ \frac{k}{k+1},\frac{m-k}{m-k+1}\}.$$  En l'infini, un changement 
sph\érique de coordonn\ées ram\ène l'\étude \à la convergence de 
\begin{eqnarray}
&&\int_{a>0}^{+\infty}\frac{
r^{-\alpha k^{2}/(k+1)}r^{-\alpha(m+1-k)}r^{m-1}} {r^{(1-\alpha)(m+1-k)+i}
r^{(1-\alpha)k-i}}dr=\int_{a>0}^{+\infty}
r^{-\alpha k^{2}/(k+1)}r^{\alpha k-2}dr;\nonumber
\end{eqnarray} 
intégrale convergente pour $\alpha < \frac{k+1}{k}$. Cette derni\ère condition 
\étant toujours v\érifi\ée, le corollaire \ref{coro2} est alors \établi.




{\subsection{Preuve du th\éor\ème \ref{th3}.}} Comme pour l'espace $X$, 
l'invariance par le
groupe $G$, des fonctions $\varphi
([z_{0},..,z_{k},z_{k+1},..,z_{m}],[\zeta_{0},..,\zeta_{k}],
[\zeta_{k+1}',..,\zeta_{m}'])$
(o\`{u} $(z_{0},..,z_{k})$ et $(\zeta_{0},..,\zeta_{k})$ ainsi que 
$(z_{k+1},..,z_{m})$ et $(\zeta_{k+1}',..,\zeta_{m}')$ sont
colin\éaires), nous permettra de les consid\érer,
dans le lemme \ref{lem31} comme des fonctions $$\varphi
([1,x_{1},...,x_{m}],[1,x_{1},...,x_{k}],[x_{k+1},...,x_{m}])$$ 
des variables r\'eelles
$ x_{i} =\mid z_{i} \mid > 0$, $i\in\{1,..,m\}$, puis, dans le lemme
\ref{lem32} comme des fonctions $$\varphi
([x_{0},...,x_{k},1,x_{k+2},..,x_{m}],[x_{0},...,x_{k}],[1,x_{k+2},...,x_{m}])$$ des
variables r\'eelles $ x_{i} =\mid z_{i} \mid >0$.
\begin{lem}\label{lem31}
Etant donn\ée une fonction $\varphi\in C^{\infty}(Y)$,
$\hat{g}$-admissible, G-invariante. Si $x_{i}=\mid z_{i}\mid >0$ pour
tout $i$:
\begin{eqnarray}\label{eq31}
&&(\varphi -\hat{\psi}
)([1,x_{1},..,x_{m}],[1,x_{1},..,x_{k}],[x_{k+1},..,x_{m}])\nonumber\\ 
&&\geq (\varphi
-\hat{\psi}
)([1,x_{1},..,x_{k};\zeta^{[m-k]}],[1,x_{1},..,x_{m}],[1^{[m-k]}]),
\end{eqnarray}
o\`u $\zeta^{[m-k]}=(\zeta ,..,\zeta)\in \complex^{m-k}$ et
$\zeta=(x_{k+1}...x_{m})^{1/(m-k)}$.
\end{lem}
{\bf Preuve.} A l'instar du lemme \ref{lem21}, nous proc\éderons par r\écurrence.
Supposons que pour $k+1\leq p < m$ et pour tout $(x_{1},..,x_{m})
\in\reel^{m}$ avec $x_{i}>0$ on ait
\begin{eqnarray}\label{eq32}
&&(\varphi -\hat{\psi})([1,x_{1},..,x_{m}],[1,x_{1},..,x_{k}],[x_{k+1},..,x_{m}])\geq \nonumber \\
&& (\varphi -\hat{\psi} )([1,x_{1},..,x_{k};(x_{k+1}...x_{p})^{\frac{1}{p-k}},..,
(x_{k+1}...x_{p})^{\frac{1}{p-k}},x_{p+1},..x_{m}],\nonumber\\
&&[1,x_{1},..,x_{k}],[(x_{k+1}...x_{p})^{\frac{1}{p-k}},..,
(x_{k+1}...x_{p})^{\frac{1}{p-k}},x_{p+1},..x_{m}]).
\end{eqnarray}
Cette propri\ét\é est claire pour $p=k+1$. Si 
l'in\égalit\é (\ref{eq32}) n'\était pas satisfaite au rang
$p+1$, il existerait alors un point
$(x_{1}^{0},..,x_{m}^{0})\in\reel^{m}$ avec $x_{i}^{0}>0$ pour
tout $i$, tel que
\begin{eqnarray}\label{eq33}
&&(\varphi -\hat{\psi} )([1,x_{1}^{0},..,x_{m}^{0}],[1,x_{1}^{0},..,x_{k}^{0}],[x_{k+1}^{0},..,x_{m}^{0}])< \nonumber \\
&&(\varphi -\hat{\psi})([1,x_{1}^{0},..,x_{k}^{0};(x_{k+1}^{0}...x_{p+1}^{0})^
{\frac{1}{p+1-k}},..,(x_{k+1}^{0}...x_{p+1}^{0})^{\frac{1}{p+1-k}},x_{p+2}^{0},..,x_{m}^{0}],\nonumber\\
&&[1,x_{1}^{0},..,x_{k}^{0}],[(x_{k+1}^{0}...x_{p+1}^{0})^
{\frac{1}{p+1-k}},..,(x_{k+1}^{0}...x_{p+1}^{0})^{\frac{1}{p+1-k}},x_{p+2}^{0},..,x_{m}^{0}]).
\end{eqnarray}

En utilisant la continuit\é de $(\varphi -\hat{\psi})$, on peut
supposer, quitte \`{a} en modifier l\ég\érement les
coordonn\ées, que le point $$([1,x_{1}^{0},..,x_{m}^{0}],[1,x_{1}^{0},..,x_{k}^{0}],
[x_{k+1}^{0},..,x_{m}^{0}])$$ de
l'in\égalit\é (\ref{eq33}), v\érifie
$$(x_{1}^{0}..x_{k}^{0})^{1/(k+1)}\neq
(x_{k+1}^{0}..x_{m}^{0})^{1/(m-k)}.$$ En utilisant la $G$-invariance de
$\varphi$, on peut supposer que $x_{k+1}^{0}\leq ...\leq
x_{m}^{0}$. D'autre part, en tenant encore compte de la $G$
invariance de $\varphi$ et de l'hypoth\èse de r\écurrence
(\ref{eq32}) en les points
\begin{eqnarray}
&&([1,x_{1}^{0},..,x_{k}^{0};x_{k+1}^{0},x_{k+2}^{0},..,x_{p}^{0},
x_{p+1}^{0},x_{p+2}^{0},..,x_{m}^{0}],\nonumber\\
&&[1,x_{1}^{0},..,x_{k}^{0}],
[x_{k+1}^{0},x_{k+2}^{0},..,x_{p}^{0},
x_{p+1}^{0},x_{p+2}^{0},..,x_{m}^{0}])\nonumber
\end{eqnarray} 
et 
\begin{eqnarray}
&&([1,x_{1}^{0},..,x_{k}^{0};x_{k+2}^{0},
x_{k+3}^{0},..,x_{p}^{0},x_{p+1}^{0},x_{k+1}^{0},x_{p+2}^{0},..,x_{m}^{0}],\nonumber\\
&&[1,x_{1}^{0},..,x_{k}^{0}],
[x_{k+2}^{0},
x_{k+3}^{0},..,x_{p}^{0},x_{p+1}^{0},x_{k+1}^{0},x_{p+2}^{0},..,x_{m}^{0}]),\nonumber
\end{eqnarray}
de $Y$, on peut \écrire
\begin{eqnarray}\label{eq34}
&&(\varphi-\hat{\psi} )([1,x_{1}^{0},..,x_{m}^{0}],[1,x_{1}^{0},..,x_{k}^{0}],
[x_{k+1}^{0},..,x_{m}^{0}])\geq\nonumber\\
&& (\varphi-\hat{\psi}
)([1,x_{1}^{0},..,x_{k}^{0};(x_{k+1}^{0}...x_{p}^{0})
^{\frac{1}{p-k}},..,(x_{k+1}^{0}...x_{p}^{0})^
{\frac{1}{p-k}},x_{p+1}^{0},x_{p+2}^{0},..,x_{m}^{0}],\nonumber\\
&&[1,x_{1}^{0},..,x_{k}^{0}],[(x_{k+1}^{0}...x_{p}^{0})
^{\frac{1}{p-k}},..,(x_{k+1}^{0}...x_{p}^{0})^
{\frac{1}{p-k}},x_{p+1}^{0},x_{p+2}^{0},..,x_{m}^{0}])
\end{eqnarray}
et
\begin{eqnarray}\label{eq35}
&&(\varphi-\hat{\psi})([1,x_{1}^{0},..,x_{k}^{0};x_{k+2}^{0},..,x_{p+1}^{0},x_{k+1}^{0},
x_{p+2}^{0},..,x_{m}^{0}],[1,x_{1}^{0},..,x_{k}^{0}],\nonumber\\
&&[x_{k+2}^{0},..,x_{p+1}^{0},x_{k+1}^{0},
x_{p+2}^{0},..,x_{m}^{0}])\geq \nonumber \\
&& (\varphi-\hat{\psi})([1,x_{1}^{0},..,x_{k}^{0};(x_{k+2}^{0}..x_{p+1}^{0})
^{\frac{1}{p-k}},..,(x_{k+2}^{0}..x_{p+1}^{0})^{\frac{1}{p-k}},x_{k+1}^{0},
x_{p+2}^{0},..,x_{m}^{0}],\nonumber\\
&&[1,x_{1}^{0},..,x_{k}^{0}],[(x_{k+2}^{0}..x_{p+1}^{0})
^{\frac{1}{p-k}},..,(x_{k+2}^{0}..x_{p+1}^{0})^{\frac{1}{p-k}},x_{k+1}^{0},
x_{p+2}^{0},..,x_{m}^{0}]).
\end{eqnarray}
Consid\érons maintenant la courbe $C$ d'\équation
$$t^{p-k}x=x_{k+1}^{0}...x_{p+1}^{0}$$ dans le plan r\éel
$$\{[1,x_{1}^{0},..,x_{k}^{0};t,..,t,x,x_{p+2}^{0},..,x_{m}^{0}],[1,x_{1}^{0},..,x_{k}^{0}],
[t,..,t,x,x_{p+2}^{0},..,x_{m}^{0}]\}$$
param\étr\é par les variables $t$ et $x$. Les points
\begin{eqnarray}
&&\hat{P}_{1}=([1,x_{1}^{0},..,x_{k}^{0};(x_{k+1}^{0}...x_{p}^{0})^
{\frac{1}{p-k}},..,(x_{k+1}^{0}...x_{p}^{0})^{\frac{1}{p-k}},
x_{p+1}^{0},x_{p+2}^{0},..,x_{m}^{0}],\nonumber\\
&&[1,x_{1}^{0},..,x_{k}^{0}],[(x_{k+1}^{0}...x_{p}^{0})^
{\frac{1}{p-k}},..,(x_{k+1}^{0}...x_{p}^{0})^{\frac{1}{p-k}},
x_{p+1}^{0},x_{p+2}^{0},..,x_{m}^{0}])\nonumber
\end{eqnarray} 
et
\begin{eqnarray}
&&\hat{P}_{2}=([1,x_{1}^{0},..,x_{k}^{0};(x_{k+2}^{0}...x_{p+1}^{0})
^{\frac{1}{p-k}},..,(x_{k+2}^{0}...x_{p+1}^{0})
^{\frac{1}{p-k}},x_{k+1}^{0},x_{p+2}^{0},..,x_{m}^{0}],\nonumber\\ 
&&[1,x_{1}^{0},..,x_{k}^{0}],[(x_{k+2}^{0}...x_{p+1}^{0})
^{\frac{1}{p-k}},..,(x_{k+2}^{0}...x_{p+1}^{0})
^{\frac{1}{p-k}},x_{k+1}^{0},x_{p+2}^{0},..,x_{m}^{0}])\nonumber
\end{eqnarray}
appartiennent \`{a} la courbe $C$. Notons que les r\éels
$x_{i}^{0}$ pour $k+1\leq i\leq p+1$ ne sont pas tous \égaux, 
sinon (\ref{eq33}) deviendrait une \égalit\é.\\ Par suite,
sachant que l'on a choisi $x_{k+1}^{0}\leq ...\leq x_{p+1}^{0}$,
les points distincts $\hat{P}_{1}$ et $\hat{P}_{2}$ se trouvent strictement de
part et d'autre de la diagonale $t=x$ du plan pr\éc\édent.
\\Or la courbe $C$ coupe cette diagonale en le point

\begin{eqnarray}
&&\hat{P}_{3}=([1,x_{1}^{0},..,x_{k}^{0};
(x_{k+1}^{0}...x_{p+1}^{0})^{\frac{1}{p+1-k}},..,
(x_{k+1}^{0}...x_{p+1}^{0})^{\frac{1}{p+1-k}},x_{p+2}^{0},..,x_{m}^{0}],\nonumber\\
&&[1,x_{1}^{0},..,x_{k}^{0}],[(x_{k+1}^{0}...x_{p+1}^{0})^{\frac{1}{p+1-k}},..,
(x_{k+1}^{0}...x_{p+1}^{0})^{\frac{1}{p+1-k}},x_{p+2}^{0},..,x_{m}^{0}])\nonumber
\end{eqnarray}
qui intervient dans l'in\égalit\é (\ref{eq33}). D'autre part,
en utilisant les relations (\ref{eq33}), (\ref{eq34}) et (\ref{eq35})
on obtient :$$ (\varphi - \hat{\psi})(\hat{P}_{3})> (\varphi - \hat{\psi}
)(\hat{P}_{1})\mbox{ et }(\varphi - \hat{\psi} )(\hat{P}_{3})> (\varphi - \hat{\psi}
)(\hat{P}_{2}),$$ ce qui prouve que la fonction $(\varphi - \hat{\psi} )$
admet un maximum local sur la courbe $C$. En cons\équence, la
restriction de la fonction $G$-invariante $(\varphi - \hat{\psi} )$
\`{a} la courbe holomorphe (toujours not\'ee $C$) d'\équation
$\xi^{p-k}z=x_{k+1}^{0}...x_{p+1}^{0}$ du plan complexe
$\{([1,x_{1}^{0},..,x_{k}^{0};\xi,..,\xi,z,x_{p+2}^{0},..,x_{m}^{0}],
[1,x_{1}^{0},..,x_{k}^{0}])\}$
atteint un maximum local en un point $\hat{P}=C(\zeta)$. Posons
$$C(\zeta)=([1,C^{1}(\zeta),..,C^{m}(\zeta)],[1,C^{1}(\zeta),..,C^{k}(\zeta)]
,[C^{k+1}(\zeta),..,C^{m}(\zeta)]),$$ 
$$\dot{C}^{\lambda}(\xi)=\frac{d C^{\lambda}}{d \xi}(\xi) \mbox{ et
} \dot{C}^{\overline{\mu}}(\xi)=\overline{\dot{C}^{\mu}(\xi)}.$$
Sachant que l'on a choisi le point $([1,x_{1}^{0},..,x_{m}^{0}],[1,x_{1}^{0},..,x_{k}^{0}],
[x_{k+1}^{0},..,x_{m}^{0}])$ de
sorte que $$(x_{1}^{0}...x_{k}^{0})^{1/(k+1)}\neq
(x_{k+1}^{0}...x_{m}^{0})^{1/(m-k)},$$ l'\équation de la
courbe $C$ et les d\éfinitions de $\hat{\psi}_{1}$ et $\hat{\psi}_{2}$ montrent 
qu'en tout point de $C$ 
\begin{eqnarray} 
\hat{\psi}_{1}([1,x_{1}^{0},..,x_{k}^{0};\xi,..,\xi,z,
x_{p+2}^{0},..,x_{m}^{0}],[1,x_{1}^{0},..,x_{k}^{0}])\neq\nonumber\\ 
\hat{\psi}_{2}([1,x_{1}^{0},..,x_{k}^{0};\xi,..,\xi,z,x_{p+2}^{0},..,x_{m}^{0}],[1,x_{1}^{0},..,x_{k}^{0}]).
\end{eqnarray}
On peut alors supposer que $\hat{\psi}=\hat{\psi}_{1}$ dans un voisinage de $\hat{P}$, 
la preuve \étant identique si l'on suppose $\hat{\psi}=\hat{\psi}_{2}$ dans
ce voisinage. On a donc :

$$\frac {\partial^{2}}{\partial \xi\partial \overline{\xi}}
\{(\varphi -\hat{\psi}_{1})(C(\zeta))\} =\frac{\partial^{2}(\varphi
-\hat{\psi}_{1})}{\partial z_{\lambda}
\partial \overline{z}_{\mu} }(C(\zeta))\dot{C}^{\lambda}(\zeta)
\dot{C}^{\overline{\mu}}(\zeta)$$ est n\égatif ou nul.
Comme $-\frac{\partial^{2}\hat{\psi}_{1}}{\partial z_{\lambda}
\partial \overline{z}_{\mu}}= g_{\lambda\overline{\mu}},$
ceci exprime que la forme hermitienne de matrice:
$$(g_{\lambda\overline{\mu}}+\frac{\partial^{2}\varphi}{\partial
z_{\lambda}
\partial \overline{z}_{\mu}})_{\lambda ,\mu}=(\frac{\partial^{2}(\varphi -\hat{\psi}_{1} )}
{\partial z_{\lambda}
\partial \overline{z}_{\mu}})_{\lambda,\mu}$$
est n\égative en $\hat{P} = C(\zeta)$. On en
d\éduit une contradiction avec la $g$-admissibilit\é de
$\varphi$ en $\hat{P}$. D'o\`u l'in\'egalit\'e (\ref{eq32}) au rang $p+1$ et, par
cons\'equent, le lemme \ref{lem31}. 


\begin{lem}\label{lem32}
Etant donn\ée une fonction $\varphi\in C^{\infty}(Y)$,
$\hat{g}$-admissible, G-invariante. Si $x_{i}=\mid z_{i}\mid >0$ pour
tout $i$, on a:
\begin{eqnarray}\label{eq36}
(\varphi -\hat{\psi}
)([x_{0},x_{1},..,x_{k};1,x_{k+2},..,x_{m}],
[x_{0},x_{1},..,x_{k}],[1,x_{k+2},..,x_{m}])\geq\nonumber\\ (\varphi -\hat{\psi}
)([\eta,\eta,..,\eta;1,x_{k+2},..,x_{m}],[1^{[k+1]}],[1,x_{k+2},..,x_{m}]),
\end{eqnarray}
o\`u $\eta=(x_{0}x_{1}...x_{k})^{1/(k+1)}$.
\end{lem}
{\bf Preuve.} Comme dans le lemme \ref{lem31}, la preuve s'effectue
par r\écurrence. Supposons que pour $0\leq p < k$ et
pour tout $(x_{0},..,x_{k};x_{k+2},..,x_{m}) \in\reel^{m}$ avec
$x_{i}>0$, on ait
\begin{eqnarray}\label{eq37}
&&(\varphi -\hat{\psi} )([x_{0},..,x_{k},1,x_{k+2},..,x_{m}],[x_{0},..,x_{k}],
[1,x_{k+2},..,x_{m}])\geq
\nonumber\\
&&(\varphi -\hat{\psi}
)([(x_{0}...x_{p})^{\frac{1}{p+1}},..,
(x_{0}...x_{p})^{\frac{1}{p+1}},x_{p+1},..,x_{k};1,x_{p+2},..,x_{m}],\nonumber\\
&&[(x_{0}...x_{p})^{\frac{1}{p+1}},..,
(x_{0}...x_{p})^{\frac{1}{p+1}},x_{p+1},..,x_{k}],[1,x_{p+2},..,x_{m}]).
\end{eqnarray}
Cette hypoth\èse est v\érifi\ée pour $p=0$. Si
l'in\égalit\é (\ref{eq37}) n'\était pas satisfaite au rang
$p+1$, il existerait un point
$(x_{0}^{0},..,x_{k};x_{k+2},..,x_{m}^{0})\in\reel^{m}$ avec
$x_{i}^{0}>0$ pour tout $i$, tel que :
\begin{eqnarray}\label{eq38}
&&(\varphi -\hat{\psi}
)([x_{0}^{0},..,x_{k}^{0};1,x_{k+2}^{0},..,x_{m}^{0}],[x_{0}^{0},..,x_{k}^{0}],[
1,x_{k+2}^{0},..,x_{m}^{0}])< \nonumber\\
&& (\varphi
-\hat{\psi})([(x_{0}^{0}..x_{p+1}^{0})^{\frac{1}{p+2}},..,(x_{0}^{0}..x_{p+1}^{0})^{\frac{1}{p+2}},
x_{p+2}^{0},...,x_{k}^{0};1,x_{k+2}^{0},..,x_{m}^{0}],\nonumber\\
&&[(x_{0}^{0}..x_{p+1}^{0})^{\frac{1}{p+2}},..,(x_{0}^{0}..x_{p+1}^{0})^{\frac{1}{p+2}},
x_{p+2}^{0},...,x_{k}^{0}],[1,x_{k+2}^{0},..,x_{m}^{0}]).
\end{eqnarray}
Comme au lemme \ref{lem31}, on peut supposer que le point
$$([x_{0}^{0},..,x_{k}^{0};1,x_{k+2}^{0},..,x_{m}^{0}],[x_{0}^{0},..,x_{k}^{0}],
[1,x_{k+2}^{0},..,x_{m}^{0}])\in Y $$ v\érifie
$$(x_{0}^{0}..x_{k}^{0})^{1/(k+1)}\neq
(x_{k+2}^{0}..x_{m}^{0})^{1/(m-k)},$$ et que $x_{0}^{0}\leq
...\leq x_{p+1}^{0}$. D'autre part, en tenant compte de la $G$-invariance 
de $\varphi$ et de l'hypoth\èse de r\écurrence
(\ref{eq37}) en les points 
$$([x_{0}^{0},x_{1}^{0},..,x_{p}^{0},x_{p+1}^{0},..,x_{k}^{0};1,x_{k+2}^{0},..,x_{m}^{0}],
[x_{0}^{0},x_{1}^{0},..,x_{p}^{0},x_{p+1}^{0},..,x_{k}^{0}],[1,x_{k+2}^{0},..,x_{m}^{0}])$$
et
$$([x_{1}^{0},..,x_{p+1}^{0},x_{0}^{0},..,x_{k}^{0};1,x_{k+2}^{0},..,x_{m}^{0}],
[x_{1}^{0},..,x_{p+1}^{0},x_{0}^{0},..,x_{k}^{0}],[1,x_{k+2}^{0},..,x_{m}^{0}]),$$
on a 
\begin{eqnarray}\label{eq39}
&&(\varphi-\hat{\psi} )([x_{0}^{0},..,x_{k}^{0};1,x_{k+2}^{0},
..,x_{m}^{0}],[x_{0}^{0},..,x_{k}^{0}],[1,x_{k+2}^{0},
..,x_{m}^{0}])\geq\nonumber\\
&& (\varphi-\hat{\psi}
)([(x_{0}^{0}..x_{p}^{0})^{\frac{1}{p+1}},..,(x_{0}^{0}..x_{p}^{0})^{\frac{1}{p+1}},
x_{p+1}^{0},..,x_{k}^{0};1,x_{k+2}^{0},..,x_{m}^{0}],\nonumber\\
&&[(x_{0}^{0}..x_{p}^{0})^{\frac{1}{p+1}},..,(x_{0}^{0}..x_{p}^{0})^{\frac{1}{p+1}},
x_{p+1}^{0},..,x_{k}^{0}],[1,x_{k+2}^{0},..,x_{m}^{0}]) 
\end{eqnarray}
et
\begin{eqnarray}\label{eq310}
&&(\varphi-\hat{\psi})([x_{1}^{0},..,x_{p+1}^{0},x_{0}^{0},..,x_{k}^{0};1,x_{k+2}^{0},..,x_{m}^{0}],\nonumber\\
&&[x_{1}^{0},..,x_{p+1}^{0},x_{0}^{0},..,x_{k}^{0}],[1,x_{k+2}^{0},..,x_{m}^{0}])\geq
\nonumber \\
&&(\varphi-\hat{\psi})([(x_{1}^{0}..x_{p+1}^{0})^{\frac{1}{p+1}},..,(x_{1}^{0}..x_{p+1}^{0})^
{\frac{1}{p+1}},
x_{0}^{0},x_{p+2}^{0},..,x_{k}^{0};1,x_{k+2}^{0},..,x_{m}^{0}],\nonumber\\
&&[(x_{1}^{0}..x_{p+1}^{0})^{\frac{1}{p+1}},..,(x_{1}^{0}..x_{p+1}^{0})^
{\frac{1}{p+1}},
x_{0}^{0},x_{p+2}^{0},..,x_{k}^{0}],[1,x_{k+2}^{0},..,x_{m}^{0}]).
\end{eqnarray}
Consid\érons maintenant la courbe $C$ d'\équation
$$t^{p+1}x=x_{0}^{0}...x_{p+1}^{0}$$ du plan r\éel
$$\{([t,..,t,x,x_{p+2}^{0},..,x_{k}^{0};1,x_{k+2}^{0},..,x_{m}^{0}],[t,..,t,x,x_{p+2}^{0},..,x_{k}^{0}],
[1,x_{k+2}^{0},..,x_{m}^{0}])\}$$
param\étr\é par les variables $t$ et $x$. Les points 
\begin{eqnarray}
&&\hat{Q}_{1}=([(x_{0}^{0}..x_{p}^{0})^{\frac{1}{p+1}},..,
(x_{0}^{0}..x_{p}^{0})^{\frac{1}{p+1}},x_{p+1}^{0},x_{p+2}^{0},..,x_{k}^{0};1,
x_{k+2}^{0},..,x_{m}^{0}],\nonumber\\
&&[(x_{0}^{0}..x_{p}^{0})^{\frac{1}{p+1}},..,
(x_{0}^{0}..x_{p}^{0})^{\frac{1}{p+1}},x_{p+1}^{0},x_{p+2}^{0},..,x_{k}^{0}],[1,
x_{k+2}^{0},..,x_{m}^{0}])\nonumber
\end{eqnarray} 
et
\begin{eqnarray}
&&\hat{Q}_{2}=([(x_{1}^{0}..x_{p+1}^{0})^{\frac{1}{p+1}},..,(x_{1}^{0}..x_{p+1}^{0})^{\frac{1}{p+1}},
x_{0}^{0},x_{p+2}^{0},..,x_{k}^{0};1,x_{k+2}^{0},..,x_{m}^{0}],\nonumber\\
&&[(x_{1}^{0}..x_{p+1}^{0})^{\frac{1}{p+1}},
..,(x_{1}^{0}..x_{p+1}^{0})^{\frac{1}{p+1}},
x_{0}^{0},x_{p+2}^{0},..,x_{k}^{0}],[1,
x_{k+2}^{0},..,x_{m}^{0}])\nonumber
\end{eqnarray} 
appartiennent \`{a} la courbe $C$.\\D'autre part 
 les r\éels $x_{i}^{0}$ pour $0\leq i\leq p+1$ ne sont pas tous
\égaux, sinon (\ref{eq38}) serait une \égalit\é.\\
Par suite, sachant que l'on a choisi $x_{0}^{0}\leq ...\leq
x_{p+1}^{0}$, les points distincts $\hat{Q}_{1}$ et $\hat{Q}_{2}$ se trouvent strictement 
de part et
d'autre de la diagonale $t=x$ du plan pr\éc\édent. \\Or la
courbe $C$ coupe cette diagonale en le point 

\begin{eqnarray}
&&\hat{Q}_{3}=([(x_{0}^{0}...x_{p+1}^{0})^{\frac{1}{p+2}},..,
(x_{0}^{0}...x_{p+1}^{0})^{\frac{1}{p+2}},x_{p+2}^{0},..,x_{k}^{0};
1,x_{k+2}^{0},..,x_{m}^{0}],\nonumber\\
&&[(x_{0}^{0}...x_{p+1}^{0})^{\frac{1}{p+2}},..,
(x_{0}^{0}...x_{p+1}^{0})^{\frac{1}{p+2}},x_{p+2}^{0},..,x_{k}^{0}],[1,
x_{k+2}^{0},..,x_{m}^{0}])\nonumber
\end{eqnarray} qui intervient dans
l'in\égalit\é (\ref{eq38}). D'autre part, les relation
(\ref{eq38}), (\ref{eq39}) et (\ref{eq310}) donnent 
$$ (\varphi -
\hat{\psi})(\hat{Q}_{3})> (\varphi - \hat{\psi} )(\hat{Q}_{1})\mbox{ et }(\varphi - \hat{\psi}
)(\hat{Q}_{3})> (\varphi - \hat{\psi} )(\hat{Q}_{2}),$$ ce qui prouve que la
fonction $(\varphi - \hat{\psi} )$ admet un maximum local sur la courbe
$C$.\\
Sachant que l'on a choisi le point
$$([x_{0}^{0},..,x_{k}^{0};1,x_{k+2}^{0},..,x_{m}^{0}],[x_{0}^{0},..,x_{k}^{0}],[1,
x_{k+2}^{0},..,x_{m}^{0}])$$ de sorte que
$$(x_{0}^{0}..x_{k}^{0})^{1/(k+1)}\neq
(x_{k+2}^{0}..x_{m}^{0})^{1/(m-k)},$$
on conclut de la
m\^{e}me mani\ère qu'au lemme pr\éc\édent en consid\érant la restriction de 
$(\varphi -\hat{\psi})$ \`{a} une courbe holomorphe convenable.


A l'instar du lemme \ref{coro24}, les lemmes \ref{lem31} et \ref{lem32} permettent d'\établir le
\begin{lem}\label{coro34}
Etant donn\ée une fonction $\varphi\in C^{\infty}(Y)$,
$\hat{g}$-admissible, G-invariante, avec $x_{i}=\mid z_{i}\mid >0$ pour
tout $i$, on a :
\begin{eqnarray}\label{eq311}
&&(\varphi -\hat{\psi}
)([1,x_{1},..,x_{k};x_{k+1},x_{k+2},..,x_{m}],[1,x_{1},..,x_{k}],[x_{k+1},x_{k+2},..,x_{m}])\nonumber\\
&&\geq(\varphi -\hat{\psi} )([1^{[k+1]};\nu^{[m-k]}],[1^{[k+1]}],[1^{[m-k]}]),
\end{eqnarray}
o\`u
$\nu=\frac{(x_{k+1}...x_{m})^{1/(m-k)}}{(x_{1}...x_{k})^{1/(k+1)}}$.
\end{lem}
Montrons maintenant le 

\begin{lem}\label{lem34}
Etant donn\ée une fonction $\varphi\in C^{\infty}(Y)$,
$\hat{g}$-admissible, G-invariante, $\forall \zeta >0$, on a :
\begin{eqnarray}\label{eq312}
(\varphi -\hat{\psi} )([1^{[k+1]};\zeta^{[m-k]}],[1^{[k+1]}],[1^{[m-k]}])\geq 0,
\end{eqnarray}
\end{lem}


{\bf Preuve.} On raisonne sur la position du point $\hat{R}_{0}\in\proj_{m}\complex$
o\`{u} $\varphi$ atteint son maximum. En vertu de la $G$-invariance
de $\varphi$, on peut supposer qu'il s'\écrit sous la forme 
$$\hat{R}_{0}=([y^{0}_{0},..,y^{0}_{k};y^{0}_{k+1},..,y^{0}_{m}],
[\rho ^{0}_{0},..,\rho ^{0}_{k}],[\varrho^{0}_{k+1},..,\varrho^{0}_{m}]),$$
o\`{u} les $y^{0}_{i}$, $\rho^{0}_{i}$, $\varrho ^{0}_{i}$ sont des r\éels positifs v\érifiant 
$y^{0}_{0}\geq y^{0}_{1}\geq..\geq y^{0}_{k}$, $y^{0}_{k+1}\geq
y^{0}_{k+2}\geq..\geq y^{0}_{m}$, $\rho ^{0}_{0}\geq..\geq\rho ^{0}_{k}$ et 
$\varrho ^{0}_{k+1}\geq..\geq\varrho ^{0}_{m}$ et o\ù $(\rho ^{0}_{0},..,\rho ^{0}_{k})$ et 
$(y^{0}_{0},..,y^{0}_{k})$ sont  
parall\èles, ainsi que $(\varrho ^{0}_{k+1},..,\varrho ^{0}_{m})$ et 
$(y^{0}_{k+1},..,y^{0}_{m})$.
Trois cas se pr\ésentent : ou
bien l'un des $y^{0}_{0},..,y^{0}_{k}$ et l'un des $y^{0}_{k},..,y^{0}_{m}$ sont non nuls, ou bien tous
les $y^{0}_{0},..,y^{0}_{k}$ sont nuls, ou bien enfin tous les $y^{0}_{k},..,y^{0}_{m}$. Ces deux derniers 
cas\étant sym\étriques, ils se traitent de mani\ère similaire.
\\{\bf{Cas A}} :
l'un des $y^{0}_{0},..,y^{0}_{k}$ et l'un des $y^{0}_{k},..,y^{0}_{m}$ sont non nuls. On se
place alors dans la carte de $Y$ d\écrite par les points 
$$([z_{0},..,z_{k},z_{k+1},..,z_{m}],[\xi_{0},..,\xi_{k}],[\xi_{k+1},..,\xi_{m}])\in
\proj_{m}\complex \times\proj_{k}\complex\times\proj_{m-k-1}\complex$$ 
tels que $z_{0}\neq 0, \xi_{0}\neq 0$.
Ceci permet de rep\érer 
$\hat{R}_{0}$ par  
$$\hat{R}_{0}=([1,x^{0}_{1},..,x^{0}_{k};x^{0}_{k+1},..,x^{0}_{m}],[1,x^{0}_{1},..,x^{0}_{k}],[
x^{0}_{k+1},..,x^{0}_{m}]),$$
o\`{u} les r\éels positifs $x^{0}_{i}$ v\érifient :
$1\geq x^{0}_{1}\geq..\geq x^{0}_{k}$ et $x^{0}_{k+1}\geq ...\geq
x^{0}_{m}$.
Raisonnons par l'absurde, et supposons qu'il existe un point
$$\hat{R}_{1}=([1^{[k+1]};\zeta_{0}^{[m-k]}],[1^{[k+1]}],[1^{[m-k]}])$$ tel que
l'on ait $\zeta_{0}>0$ et 
\begin{equation}\label{p31}
(\varphi -\hat{\psi} )(\hat{R}_{1})<0.
\end{equation}
On envisage alors les deux sous-cas suivants : $x^{0}_{k+1}<\zeta_{0}$ 
puis $x^{0}_{k+1}\geq \zeta_{0}$.
\begin{itemize}
    \item \underline{$x^{0}_{k+1}<\zeta_{0}$}.
\end{itemize}
On introduit alors la fonction auxiliaire
$$\hat{\psi}_{0}=\log\frac{\mid z_{0}\mid^{4}\mid \xi_{0}\mid^{2k}\mid \xi_{k+1}\mid^{2(m-k-1)}}{
    (\mid z_{0}\mid^{2}+...+\mid z_{m}\mid^{2})^{2}(\mid \xi_{0}\mid^{2}+...+\mid \xi_{k}\mid^{2})^{k}
(\mid \xi_{k+1}\mid^{2}+...+\mid \xi_{m}\mid^{2})^{m-k-1}}.$$
D'une part, puisque $\varphi\leq 0$,
\begin{equation}\label{p30}
(\varphi -\hat{\psi}_{0})([1,0^{[m]}],[1,0^{[k]}],[1,0^{[m-k-1]}])=\varphi 
([1,0^{[m]}],[1,0^{[k]}],[1,0^{[m-k-1]}])\leq 0.
\end{equation}
De plus, sachant que $\varphi (\hat{R}_{0})=0$ et $\hat{\psi}_{0}\leq 0$, on a
\begin{equation}\label{p33}
(\varphi -\hat{\psi}_{0})(\hat{R}_{0})\geq 0.
\end{equation}
Si $\hat{R}_{0}\neq ([1,0^{[m]}],[1,0^{[k]}],[1,0^{[m-k-1]}])$, $\hat{\psi}_{0}(\hat{R}_{0})< 0$ 
et l'in\égalit\é (\ref{p33}) est
alors stricte. Si $\hat{R}_{0} = ([1,0^{[m]}],[1,0^{[k]}],[1,0^{[m-k-1]}])$, 
quitte \`{a} se placer en un point $\hat{R}$ arbitrairement voisin de
$\hat{R}_{0}$, on peut supposer $(\varphi -\hat{\psi}_{0})(\hat{R})> 0$. En effet, si dans
un voisinage de $\hat{R}_{0}$ on avait $(\varphi -\hat{\psi}_{0})\leq 0$, comme 
$(\varphi -\hat{\psi}_{0})(\hat{R}_{0})= 0$, $(\varphi -\hat{\psi}_{0})$ admettrait
alors un maximun local en $\hat{R}_{0}$, ce qui mettrait en
d\éfaut l'admissibilit\é de $\varphi$ en ce point, sachant que 
$$\partial_{\lambda\overline{\mu}}(\varphi -\hat{\psi}_{0})(\hat{R}_{0})=
(g_{\lambda\overline{\mu}}+\partial_{\lambda\overline{\mu}}\varphi )(\hat{R}_{0}).$$ Dans tous les cas, 
on peut
donc affirmer qu'il existe un point $$\hat{R}'_{0}=([1,{a}_{1},..,{a}_{m}],[1,{a}_{1},..,{a}_{k}]
,[{a}_{k+1},..,{a}_{m}])$$
v\érifiant
\begin{equation}\label{p33'}
(\varphi -\hat{\psi}_{0})(\hat{R}'_{0})> 0.
\end{equation}
Par continuit\é et $G$-invariance de $\varphi$, on peut supposer $1>a_{1}>...> a_{k}>0$ et 
$\zeta_{0}> a_{k+1}>...> a_{m}>0$.
D'autre part, l'in\égalit\é (\ref{p31}) jointe aux d\'efinitions de $\hat{R}_1$, $\hat{\psi}_0$, $\hat{\psi}_1$ 
et $\hat{\psi}=\inf (\hat{\psi}_1 , \hat{\psi}_2)$ implique
\begin{equation}\label{p31'}
(\varphi -\hat{\psi}_{0})(\hat{R}_{1})=(\varphi -\hat{\psi}_{1})(\hat{R}_{1})\leq(\varphi
-\hat{\psi})(\hat{R}_{1})< 0.
\end{equation}

La courbe :
\begin{eqnarray}
&[0,1]\ni t\rightarrow ([1,t,t^{(\ln a_{2})/(\ln a_{1})},..,t^{(\ln a_{k})/(\ln a_{1})};
\zeta_{0}t^{\frac{\ln(a_{k+1}/\zeta_{0})}{\ln
a_{1}}},..,\zeta_{0}t^{\frac{\ln(a_{m}/\zeta_{0})}{\ln
a_{1}}}],&\nonumber\\
&[1,t,t^{(\ln a_{2})/(\ln a_{1})},..,t^{(\ln a_{k})/(\ln a_{1})}],[1,t^{\frac{\ln(a_{k+2}/a_{k+1})}{\ln
a_{1}}},..,t^{\frac{\ln(a_{m}/a_{k+1})}{\ln
a_{1}}}])&
\nonumber
\end{eqnarray} 
passe par $([1,0^{[m]}],[1,0^{[k]}],[1,0^{[m-k-1]}])$ en $t=0$ puis par ${\hat{R}'}_{0}$
en $t=a_{1}$ et enfin par le point $\hat{R}_{1}$ en $t=1$, valeurs
en lesquelles, d'apr\ès (\ref{p30}), (\ref{p33'}) et (\ref{p31'}),
$(\varphi-\hat{\psi}_{0})$ est respectivement n\égative, positive
puis \`{a} nouveau n\égative. L'invariance de cette fonction
par l'action des $\exp(i\theta)$, permet donc de d\éduire que
$(\varphi-\hat{\psi}_{0})$ atteint un maximum sur la courbe holomorphe,
complexifi\ée de la courbe d\écrite plus haut, ce qui
contredit encore une fois l'admissibilit\é de $\varphi$.
\begin{itemize}
\item \underline{$x^{0}_{k+1}\geq\zeta_{0}$}.
\end{itemize}
D\ésignons dans ce cas par $p\in \{1,..,m-k\}$ l'entier pour
lequel on a $$x^{0}_{k+1}\geq ...\geq x^{0}_{k+p}>\zeta_{0} \mbox{
et } \zeta_{0}\geq x^{0}_{k+p+1}\geq ...\geq x^{0}_{m},$$ et
consid\érons la fonction auxiliaire

$$\hat{\psi}_{k+1}=\log\frac{\mid z_{k+1}\mid^{4}\mid\xi_{0}\mid^{2k}\mid\xi_{k+1}\mid^{2(m-k-1)}}{
(\mid z_{0}\mid^{2}+...+\mid z_{m}\mid^{2})^{2}(\mid \xi_{0}\mid^{2}+...+\mid \xi_{k}\mid^{2})^{k}
(\mid \xi_{k+1}\mid^{2}+...+\mid \xi_{m}\mid^{2})^{m-k-1}}.$$
On a
\begin{equation}\label{p33''}
(\varphi -\hat{\psi}_{k+1})(\hat{R}_{0})> 0.
\end{equation}
La fonction $(\varphi -\hat{\psi}_{k+1})$ \étant continue, quitte
\`{a} se placer en un point voisin de $\hat{R}_{0}$, on peut supposer
tous les $x_{i}^{0}$ non nuls. Posons alors:
$$\alpha_{2}=\frac{\ln x_{2}^{0}}{\ln x_{1}^{0}},..,\alpha_{k}=\frac{\ln x_{k}^{0}}{\ln x_{1}^{0}};
\alpha_{k+1}=\frac{\ln(x_{k+1}^{0}/\zeta_{0})}{\ln x_{1}^{0}},..,
\alpha_{m}=\frac{\ln(x_{m}^{0}/\zeta_{0})}{\ln x_{1}^{0}}.$$
Sachant que l'on a $1\geq x_{1}^{0}\geq..\geq x_{k}^{0}$;
$x_{k+1}^{0}\geq..\geq x_{k+p}^{0}\geq\zeta_{0}$ et $\zeta_{0}\geq
x_{k+p+1}^{0}\geq..\geq x_{m}^{0}$, on en d\éduit que 
$\alpha_{2},..,\alpha_{k}\geq 0$, $\alpha_{k+1}\leq ...\leq
\alpha_{p+k}\leq 0$ et $\alpha_{p+k+1},..,\alpha_{m}\geq 0$, d'o\`{u}, en notant 
$$\hat{R}_{\varepsilon}=([1,\varepsilon,\varepsilon^{\alpha_{2}},..,\varepsilon^{\alpha_{k}};
\zeta_{0}\varepsilon^{\alpha_{k+1}},..,\zeta_{0}\varepsilon^{\alpha_{m}}],
[1,\varepsilon,\varepsilon^{\alpha_{2}},..,\varepsilon^{\alpha_{k}}],
[\varepsilon^{\alpha_{k+1}},..,\varepsilon^{\alpha_{m}}]),$$ on a 
\begin{eqnarray*}
\lim_{\varepsilon\rightarrow
0}\hat{\psi}_{k+1}(\hat{R}_{\varepsilon }) 
&=& \lim_{\varepsilon\rightarrow
0}\{\ln\frac{\zeta_{0}^{4}\varepsilon^{4\alpha_{k+1}}}
{[1+\varepsilon^{2}+\varepsilon^{2\alpha_{2}}+..+\varepsilon^{2\alpha_{k}}+
\zeta_{0}^{2}(\varepsilon^{2\alpha_{k+1}}+..+\varepsilon^{2\alpha_{m}})]^{2}}  \\
&&\times\frac{\varepsilon^{2(m-k-1)\alpha_{k+1}}}
{(\varepsilon^{2\alpha_{k+1}}+..+\varepsilon^{2\alpha_{m}})^{m-k-1}}\}\\
&=&\ln\lim_{t\rightarrow
\infty}\frac{t^{-2\alpha_{k+1}(m+1-k)}}
{(t^{-2\alpha_{k+1}}+t^{-2\alpha_{k+2}}+..+t^{-2\alpha_{p}})^{(m+1-k)}}
= \ln 1=0,
\end{eqnarray*} 
$-\alpha_{k+1}$ \étant la plus grande puissance intervenant au
d\énominateur. Sachant que $\varphi ([\hat{R}_{\varepsilon}])\leq 0$ et compte tenu de (\ref{p33''}) il
existe $\varepsilon_{0}$ tel que l'on ait
\begin{equation}\label{p34}
(\varphi
-\hat{\psi}_{k+1})(\hat{R}_{\varepsilon_{0} })\leq -\hat{\psi}_{k+1}(\hat{R}_{\varepsilon_{0} })<
(\varphi -\hat{\psi}_{k+1})(\hat{R}_{0}).
\end{equation}
D'autre part, l'in\'egalit\'e (\ref{p31}), jointe aux d\'efinitions de $\hat{R}_1$, $\hat{\psi}_{k+1}$, $\hat{\psi}_{2}$ 
et $\hat{\psi} = \inf (\hat{\psi}_{1},\hat{\psi}_{2})$ donne~: 
\begin{equation}\label{p35}
(\varphi
-\hat{\psi}_{k+1})(\hat{R}_{1})= (\varphi-\hat{\psi}_{2})(\hat{R}_{1})\leq
(\varphi -\hat{\psi})(\hat{R}_{1})<0.
\end{equation}
La courbe
$$[\varepsilon_{0},1]\ni t\rightarrow ([1,t,t^{\alpha_{2}},..,t^{\alpha_{k}},
\zeta_{0}t^{\alpha_{k+1}},..,\zeta_{0}t^{\alpha_{m}}],
[1,t,t^{\alpha_{2}},..,t^{\alpha_{k}}],
[t^{\alpha_{k+1}},..,t^{\alpha_{m}}]),$$ passe par
$\hat{R}_{\varepsilon_{0}}$
en $t=\varepsilon_{0}$ puis par $\hat{R}_{0}$ en $t=x_{1}^{0}$ et enfin
par $\hat{R}_{1}$ en $t=1$, ce qui, en vertu de 
(\ref{p34}), (\ref{p33''}) et (\ref{p35}) prouve l'existence d'un maximum local pour la 
fonction $(\varphi-\hat{\psi}_{k+1})$ sur la
courbe pr\écit\ée. Ceci contredit, \`{a} l'instar du cas
pr\éc\édent, l'hypoth\èse d'admissibilit\é de la
fonction $\varphi$.

{\bf{Cas B :}}  $y^{0}_{0}=...=y^{0}_{k}=0$. On peut alors se
placer dans la carte de $Y$ d\écrite par les points 
$$([z_{0},..,z_{k};z_{k+1},..,z_{m}],[\xi_{0},..,\xi_{k}],[\xi_{k+1},..,\xi_{m}])\in
\proj_{m}\complex \times\proj_{k}\complex\times\proj_{m-k-1}\complex$$ 
tels que $z_{k+1}\neq 0$ et $\xi_{0}\neq 0$, de sorte que
le point $\hat{R}_{0}$ o\`u $\varphi$ atteint son maximum \'egal \`a z\'ero puisse 
s'\'ecrire sous la forme 
$$\hat{R}_{0}=([0,0,..,0;1,x_{k+2}^{0},..,x_{m}^{0}],[1,u_{1}^{0},..,u_{k}^{0}],
[1,x_{k+2}^{0},..,x_{m}^{0}]).$$
On peut aussi supposer, en utilisant la $G$-invariance de $\varphi$, que 
$1\geq x_{k+2}^{0}\geq..\geq x_{m}^{0}$ et $1\geq u_{1}^{0}\geq..\geq u_{k}^{0}$. 
On montrera une version \'equivalente du lemme \ref{lem34}, \`a savoir que
\begin{equation}\label{equ30}
(\varphi -\hat{\psi} )([\zeta^{[k+1]};1^{[m-k]}],[1^{[k+1]}],[1^{[m-k]}])\geq 0
\end{equation}
pour tout $\zeta >0$. 

Raisonnons par l'absurde, et supposons qu'il existe un point
$$\hat{R}_{k+1}=([\zeta_{0}^{[k+1]};1^{[m-k]}],[1^{[k+1]}],[1^{[m-k]}])$$ de $Y$ tel que
l'on ait $\zeta_{0}>0$ et 
\begin{equation}\label{equ31}
(\varphi -\hat{\psi} )(\hat{R}_{k+1})<0.
\end{equation}
On consid\`ere alors la fonction auxiliaire
$\hat{\psi}_{k+1}$.
Sachant que $\varphi\leq 0$, on a 
\begin{eqnarray}\label{equ32}
&&(\varphi -\hat{\psi}_{k+1})([0^{[k+1]};1,0,..,0],[1,0^{[k]}],[1,0^{[m-k-1]}])=\nonumber\\
&&\varphi ([0^{[k+1]};1,0,..,0],[1,0^{[k]}],[1,0^{[m-k-1]}])\leq 0.
\end{eqnarray}
D'autre part, sachant que $\varphi (\hat{R}_{0})=0$ et $\hat{\psi}_{k+1}\leq 0$, on a

\begin{equation}
(\varphi -\hat{\psi}_{k+1})(\hat{R}_{0})=-\hat{\psi}_{k+1}(\hat{R}_{0})\geq 0,
\end{equation}
in\'egalit\'e stricte d\`es que $$\hat{R}_{0}\neq ([0^{[k+1]};1,0,..,0],[1,0^{[k]}],[1,0^{[m-k-1]}]).$$ Si 
$\hat{R}_{0} = ([0^{[k+1]};1,0,..,0],[1,0^{[k]}],[1,0^{[m-k-1]}])$, quitte \`a se placer en un point arbitrairement voisin de 
$\hat{R}_{0}$, on peut supposer la derni\`ere in\'egalit\'e stricte. En effet, si dans un voisinage de $\hat{R}_{0}$, on avait 
$\varphi -\hat{\psi}_{k+1}\leq 0$, alors $\varphi -\hat{\psi}_{k+1}$ admettrait un maximum local en  
$\hat{R}_{0}$, ce qui contredirait l'admissibilit\'e de $\varphi$ en $\hat{R}_{0}$. A l'instar du cas A, 
il existe donc un point  $$\hat{R}'_{0}=([{c}_{0},..,{c}_{k};1,{c}_{k+2},..,{c}_{m}],[{c}_{0},..,{c}_{k}]
,[1,{c}_{k+2},..,{c}_{m}])$$
v\érifiant
\begin{equation}\label{equ33}
(\varphi -\hat{\psi}_{k+1})(\hat{R}'_{0})> 0.
\end{equation}
Par continuit\é et $G$-invariance de $\varphi$, on peut supposer 
$\zeta_{0}> c_{0}>...> c_{k}>0$ et $1>c_{k+2}>...> c_{m}>0$. 
D'autre part, l'in\égalit\é (\ref{equ31}) jointe aux d\'efinitions de 
$\hat{R}_{k+1}$, $\hat{\psi}_{k+1}$, $\hat{\psi}_2$ 
et $\hat{\psi}=\inf (\hat{\psi}_1 , \hat{\psi}_2)$ implique
\begin{equation}\label{equ34}
(\varphi -\hat{\psi}_{k+1})(\hat{R}_{k+1})=(\varphi -\hat{\psi}_{2})(\hat{R}_{k+1})\leq(\varphi
-\hat{\psi})(\hat{R}_{k+1})< 0.
\end{equation}

La courbe de $Y$~:
\begin{eqnarray}
&[0,1]\ni t\rightarrow ([\zeta_{0}t^{\frac{\ln(c_{0}/\zeta_{0})}{\ln
c_{k+2}}},..,\zeta_{0}t^{\frac{\ln(c_{k}/\zeta_{0})}{\ln
c_{k+2}}};1,t,t^{(\ln c_{k+3})/(\ln c_{k+2})},..,t^{(\ln c_{m})/(\ln c_{k+2})}],&\nonumber\\
&[1,t^{\frac{\ln(c_{1}/c_{0})}{\ln
c_{k+2}}},..,t^{\frac{\ln(c_{k}/c_{0})}{\ln
c_{k+2}}}],[1,t,t^{(\ln c_{k+3})/(\ln c_{k+2})},..,t^{(\ln c_{m})/(\ln c_{k+2})}])&\nonumber
\end{eqnarray} 
passe par $([0^{[k+1]},1,0,..,0],[1,0^{[k]}])$ en $t=0$ puis par $\hat{R}_{0}$
en $t=c_{k+2}$ et enfin par le point $\hat{R}_{k+1}$ en $t=1$, valeurs
en lesquelles, d'apr\ès (\ref{equ32}), (\ref{equ33}) et (\ref{equ34}),
$(\varphi-\hat{\psi}_{k+1})$ est respectivement n\égative, positive
puis n\égative. L'invariance de cette fonction
par l'action des $\exp(i\theta)$, permet donc de d\éduire que
$(\varphi-\hat{\psi}_{0})$ atteint un maximum sur la courbe holomorphe,
d\'eduite de la courbe d\écrite plus haut, ce qui
contredit l'admissibilit\é de $\varphi$, d'o\`u (\ref{equ30}).

{\bf{Cas C :}} Les $y^{0}_{0},..,y^{0}_{k}$ sont tous nuls. On peut alors se
placer dans la carte de $Y$ d\écrite par les points 
$$([z_{0},..,z_{k};z_{k+1},..,z_{m}],[\xi_{0},..,\xi_{k}],[\xi_{k+1},..,\xi_{m}])\in
\proj_{m}\complex \times\proj_{k}\complex\times\proj_{m-k-1}\complex$$ 
tels que $z_{0}\neq 0$ et $\xi_{k+1}\neq 0$, de sorte que
le point $\hat{R}_{0}$ o\`u $\varphi$ atteint son maximum \'egal \`a z\'ero puisse 
s'\'ecrire sous la forme 
$$\hat{R}_{0}=([1,x_{1}^{0},..,x_{k}^{0};0,0,..,0],[1,x_{1}^{0},..,x_{k}^{0}],
[1,u_{k+2}^{0},..,u_{m}^{0}]).$$
On peut aussi supposer, en utilisant la $G$-invariance de $\varphi$, que 
$1\geq x_{1}^{0}\geq..\geq x_{k}^{0}$ et $1\geq u_{k+2}^{0}\geq..\geq u_{m}^{0}$.

Raisonnons par l'absurde, et supposons qu'il existe un point que l'on notera encore 
$$\hat{R}_{k+1}=([1^{[k+1]};\zeta_{0}^{[m-k]}],[1^{[k+1]}],[1^{[m-k]}])$$ de $Y$ tel que
l'on ait $\zeta_{0}>0$ et 
\begin{equation}\label{equ31'}
(\varphi -\hat{\psi} )(\hat{R}_{k+1})<0.
\end{equation}
On consid\`ere alors la fonction auxiliaire
$\hat{\psi}_{0}$ d\éfinie plus haut.
Sachant que $\varphi\leq 0$, on a 
\begin{eqnarray}\label{equ32'}
&&(\varphi -\hat{\psi}_{0})([1,0,..,0;0^{[m-k]}],[1,0^{[k]}],[1,0^{[m-k-1]}])=\nonumber\\
&&\varphi ([1,0^{[k]};0^{[m-k]}],[1,0^{[k]}],[1,0^{[m-k-1]}])\leq 0.
\end{eqnarray}
D'autre part, sachant que $\varphi (\hat{R}_{0})=0$ et $\hat{\psi}_{k+1}\leq 0$, on a

\begin{equation}
(\varphi -\hat{\psi}_{0})(\hat{R}_{0})=-\hat{\psi}_{k+1}(\hat{R}_{0})\geq 0,
\end{equation}
in\'egalit\'e stricte d\`es que $$\hat{R}_{0}\neq ([1,0^{[k]};0^{[m-k]}],[1,0^{[k]}],[1,0^{[m-k-1]}]).$$ Si 
$([1,0^{[k]};0^{[m-k]}],[1,0^{[k]}],[1,0^{[m-k-1]}])$, quitte \`a se placer en un point arbitrairement voisin de 
$\hat{R}_{0}$, on peut supposer la derni\`ere in\'egalit\'e stricte. En effet, si dans un voisinage de $\hat{R}_{0}$, on avait 
$\varphi -\hat{\psi}_{k+1}\leq 0$, alors $\varphi -\hat{\psi}_{k+1}$ admettrait un maximum local en  
$\hat{R}_{0}$, ce qui contredirait l'admissibilit\'e de $\varphi$ en $\hat{R}_{0}$.  
Il existe donc un point  $$\hat{R}'_{0}=([1,{c}_{1},..,{c}_{k};{c}_{k+1},..,{c}_{m}],[1,{c}_{1},..,{c}_{k}]
,[{c}_{k+1},..,{c}_{m}])$$
v\érifiant
\begin{equation}\label{equ33'}
(\varphi -\hat{\psi}_{0})(\hat{R}'_{0})> 0.
\end{equation}
Par continuit\é et $G$-invariance de $\varphi$, on peut supposer 
$\zeta_{0}> c_{k+1}>...> c_{m}>0$ et $1>c_{1}>...> c_{k}>0$. 
D'autre part, l'in\égalit\é (\ref{equ31'}) jointe aux d\'efinitions de 
$\hat{R}_{k+1}$, $\hat{\psi}_{0}$, $\hat{\psi}_1$ 
et $\hat{\psi}=\inf (\hat{\psi}_1 , \hat{\psi}_2)$ implique
\begin{equation}\label{equ34'}
(\varphi -\hat{\psi}_{0})(\hat{R}_{k+1})=(\varphi -\hat{\psi}_{1})(\hat{R}_{k+1})\leq(\varphi
-\hat{\psi})(\hat{R}_{k+1})< 0.
\end{equation}

La courbe de $Y$~:
\begin{eqnarray}
&[0,1]\ni t\rightarrow ([1,t,t^{(\ln c_{2})/(\ln c_{1})},..,t^{(\ln c_{k})/(\ln c_{1})};
\zeta_{0}t^{\frac{\ln(c_{k+1}/\zeta_{0})}{\ln
c_{1}}},..,\zeta_{0}t^{\frac{\ln(c_{m}/\zeta_{0})}{\ln
c_{1}}}],&\nonumber\\
&[1,t,t^{(\ln c_{2})/(\ln c_{1})},..,t^{(\ln c_{k})/(\ln c_{1})}],[1,t^{\frac{\ln(c_{k+2}/c_{k+1})}{\ln
c_{1}}},..,t^{\frac{\ln(c_{m}/c_{k+1})}{\ln
c_{1}}}])&
\nonumber
\end{eqnarray} 
passe par $([1,0,..,0],[1,0^{[k]}],[1,0^{[m-k]}])$ en $t=0$ puis par $\hat{R}_{0}$
en $t=c_{1}$ et enfin par le point $\hat{R}_{k+1}$ en $t=1$, valeurs
en lesquelles, d'apr\ès (\ref{equ32'}), (\ref{equ33'}) et (\ref{equ34'}),
$(\varphi-\hat{\psi}_{k+1})$ est respectivement n\égative, positive
puis n\égative. L'invariance de cette fonction
par l'action des $\exp(i\theta)$, permet donc de d\éduire que
$(\varphi-\hat{\psi}_{0})$ atteint un maximum sur la courbe holomorphe,
d\'eduite de la courbe d\écrite plus haut, ce qui
contredit l'admissibilit\é de $\varphi$, d'o\`u 
et le lemme \ref{lem34}.


{\subsection{Preuve du corollaire \ref{coro3}.}}

Soit $\varphi\in C^{\infty}(Y)$ une fonction
$\hat{g}$-admissible et $G$-invariante, dont le $\sup$ sur
$Y$ est nul. D'apr\`es le th\'eor\`eme \ref{th3},
on a $\varphi \geq \hat{\psi}$ et par suite, pour tout $\alpha \geq 0$,
$$\int_{Y} \exp(-\alpha\varphi ) dv \leq
\int_{Y}\exp(-\alpha\hat{\psi} )dv .$$ Afin d'obtenir les
valeurs de $\alpha$ pour lesquelles cette derni\`ere int\'egrale
converge, on estimera
$\int_{Y}\exp(-\alpha\hat{\psi}_{1} )dv$ et $\int_{Y}\exp(-\alpha\hat{\psi}_{2} )dv$ 
dans la carte
dense correspondant \à la param\étrisation 
$$([1,z_{1},..,z_{m}],[1,z_{1},..,z_{k}],[z_{k+1},..,z_{m}]).$$ Dans cette carte,
l'\'el\'ement de volume est donn\'e par (c.f. \cite{BC})~: 
$$dv=\det((\hat{g}_{\lambda \overline{\mu}}))dz_{1}\wedge d\overline{z}_{1}\wedge ...\wedge
dz_{m}\wedge d\overline{z}_{m},$$ o\ù
\begin{eqnarray}
\det((\hat{g}_{\lambda \overline{\mu}}))=2(-1)^{m}\frac{[(k+2)(1+\mid z_{1}\mid^{2}+...+\mid
z_{k}\mid^{2})+k(\mid z_{k+1}\mid^{2}+...+\mid
z_{m}\mid^{2})]^{k}} {(1+\mid z_{1}\mid^{2}+...+\mid
z_{k}\mid^{2})^{k}(\mid z_{k+1}\mid^{2}+...+\mid
z_{m}\mid^{2})^{m-k-1}}\nonumber\\
\times \frac {[(m-k-1)(1+\mid z_{1}\mid^{2}+...+\mid
z_{k}\mid^{2})+(m-k+1)(\mid z_{k+1}\mid^{2}+...+\mid
z_{m}\mid^{2})]^{m-k-1}}{(1+\mid z_{1}\mid^{2}+...+\mid
z_{m}\mid^{2})^{m+1}}.\nonumber
\end{eqnarray} 
Si $u_{p}=\mid z_{p}\mid^{2}$, en z\éro et en l'infini on a l'\équivalence~: 
$$ dv\sim \frac{Cst(>0)}{(1+u_{1}+..+u_{k})^{k}(u_{k+1}+..+u_{m})^{m-k-1}(1+u_{1}+..+u_{m})^{2}}.$$
La convergence de $\int_{Y}\exp(-\alpha\hat{\psi}_{1})dv$ est donc \équivalente \à celle de
\begin{eqnarray}
\int_{0}^{+\infty}..\int_{0}^{+\infty}\frac{
(u_{1}...u_{k})^{-\alpha(k+2)/(k+1)}(u_{k+1}...u_{m})^{-\alpha(m-k-1)/(m-k)}} {(1+
u_{1}+...+u_{m})^{2(1-\alpha)}
(1+u_{1}+..+u_{k})^{(1-\alpha)k}}\nonumber\\
\times\frac{du_{1}...du_{m}}{(u_{k+1}+..+u_{m})^{(1-\alpha)(m-k-1)}},
\end{eqnarray}
qui converge en z\éro pour $\alpha<(k+1)/(k+2)$. En l'infini en effectuant en changement de 
coordonn\ées sph\ériques on ram\ène l'\étude \à celle de 
\begin{eqnarray}
\int_{a>0}^{+\infty}\frac{
r^{-\alpha k(k+2)/(k+1)}r^{-\alpha(m-1-k)}r^{m-1}} {r^{2(1-\alpha)}r^{k(1-\alpha)}
r^{(1-\alpha)(m-k-1)}}dr,
\end{eqnarray} 
qui converge pour $\alpha < (k+1)/(k+2)$. 
De m\ême, la convergence de $\int_{Y}\exp(-\alpha\hat{\psi}_{2})dv$ est \équivalente \à celle de
\begin{eqnarray}
\int_{0}^{+\infty}..\int_{0}^{+\infty}\frac{
(u_{1}...u_{k})^{-\alpha k/(k+1)}(u_{k+1}...u_{m})^{-\alpha(m-k+1)/(m-k)}} {(1+
u_{1}+...+u_{m})^{2(1-\alpha)}
(1+u_{1}+..+u_{k})^{(1-\alpha)k}}\nonumber\\
\times\frac{du_{1}...du_{m}}{(u_{k+1}+..+u_{m})^{(1-\alpha)(m-k-1)}},
\end{eqnarray}
qui converge en z\éro pour $\alpha<1/2$. En l'infini, cela revient \à \étudier la convergence de
\begin{eqnarray}
\int_{a>0}^{+\infty}
r^{-\alpha k^{2}/(k+1)}r^{-\alpha (m-k+1)}r^{-(1-\alpha)(m+1)} r^{m-1}dr,
\end{eqnarray} 
qui converge pour $\alpha < (k+1)/k$, condition toujours v\érifi\ée, d'o\ù le corollaire \ref{coro3}.

\begin{center}
--------------------------

Universit\'e Pierre et Marie Curie, Paris, France.

e-mail: benabdes@math.jussieu.fr \\
           dridi@math.jussieu.fr

\end{center}

\end{document}